\documentclass[10pt]{article}
\usepackage{amsmath, amsthm}
\usepackage{amscd}
\usepackage{amsfonts}
\usepackage{amssymb}
\usepackage{float}
\usepackage{mathrsfs}
\usepackage{tikz}
\usepackage{color}
\definecolor{red}{rgb}{.7,0,0}
\definecolor{blue}{rgb}{0,0,1}

\usepackage[labelfont=bf,labelsep=space]{caption}

\newtheorem{theorem}{Theorem}[section]
\newtheorem{proposition}[theorem]{Proposition}
\newtheorem{lemma}[theorem]{Lemma}
\newtheorem{corollary}[theorem]{Corollary}

\theoremstyle{definition}
\newtheorem{definition}[theorem]{Definition}

\newtheorem{remark}[theorem]{Remark}

\def\proof{{\noindent\sc Proof. \quad}}
\def\sketchproof{{\noindent\sc Sketch of proof. \quad}}
\newcommand\proofof[1]{{\noindent\sc Proof of #1. \quad}}
\def\eproof{{\mbox{}\hfill\qed}\medskip}

\begin{document}

\makeatletter

\renewcommand{\bar}{\overline}

\newcommand{\x}{\times}
\newcommand{\<}{\langle}
\renewcommand{\>}{\rangle}
\newcommand{\into}{\hookrightarrow}

\renewcommand{\a}{\alpha}
\renewcommand{\b}{\beta}
\renewcommand{\d}{\delta}
\newcommand{\D}{\Delta}
\newcommand{\e}{\varepsilon}
\newcommand{\g}{\gamma}
\newcommand{\G}{\Gamma}
\renewcommand{\l}{\lambda}
\renewcommand{\L}{\Lambda}
\newcommand{\n}{\nabla}
\newcommand{\var}{\varphi}
\newcommand{\s}{\sigma}
\newcommand{\Sig}{\Sigma}
\renewcommand{\t}{\theta}
\renewcommand{\O}{\Omega}
\renewcommand{\o}{\omega}
\newcommand{\z}{\zeta}
\newcommand{\balpha}{\boldsymbol \alpha}
\newcommand{\ab}{\alpha_{\bullet}}
\def\ba{\ol{\a}}
\def\bb{\ol{\b}}
\def\bg{\ol{\gamma}}
\def\bt{\ol{\tau}}
\def\mun{\mu_{\mathrm{norm}}}
\def\sep{\mathsf{sep}}
\def\Tub{\mathrm{Tub}}
\def\hmg{\mathrm{proj}}

\newcommand{\p}{\partial}
\renewcommand{\hat}{\widehat}
\renewcommand{\bar}{\overline}
\renewcommand{\tilde}{\widetilde}


\font\eightrm=cmr8
\font\ninerm=cmr9


\def\N{\mathbb{N}}
\def\Z{\mathbb{Z}}
\def\R{\mathbb{R}}
\def\Q{\mathbb{Q}}
\def\C{\mathbb{C}}
\def\F{\mathbb{F}}
\def\IS{\mathbb{S}}
\def\proj{\mathbb{P}}
\def \Ri{\R^\infty}
\def \Zi{\Z^\infty}
\def \ZRi{\Z^\infty\x\R^\infty}
\def \SZi{\Z^\infty\x S(\R^\infty)}


\newcommand{\algoritmo}{\begin{minipage}{0.87\hsize}\linea}
\newcommand{\falgoritmo}{\linea\end{minipage}\bigskip}
\newcommand{\codigo}{\begin{minipage}{0.87\hsize}}
\newcommand{\fcodigo}{\end{minipage}\bigskip}
\newcommand{\linea}{\vspace*{-5pt}\hrule\vspace*{5pt}}
\newtheorem{algorithm}{Algorithm}
\def\espacio{\hspace*{1cm}}
\def\eespacio{\hspace*{1.5cm}}
\def\eeespacio{\hspace*{2cm}}
\newcommand{\inputalg}[1]{\linea\bf Input:\quad\rm #1\vspace*{3pt}}
\newcommand{\specalg}[1]{\bf Preconditions:\quad\rm #1}
\newcommand{\Output}[1]{\linea\bf Output:\quad\rm #1\vspace*{2pt}}
\newcommand{\postcond}[1]{\bf Postconditions:\quad\rm #1\vspace*{3pt}}
\newcommand{\bodyalg}[1]{\linea\tt #1\vspace*{3pt}}
\newcommand{\bodycode}[1]{\tt #1\vspace*{3pt}}
\newcommand{\sssection}[1]{\ \\ \noindent{\bf #1.}\quad}


\def\sg{{\rm sign}\,}
\def\sig{\overline{\rm sign}\,}
\def\rank{{\rm rank}\,}
\def\Im{{\rm Im}}
\def\dist{{\rm dist}}
\def\degree{{\rm degree}\,}
\def\grad{{\rm grad}\,}
\def\size{{\rm size}}
\def\supp{{\rm supp}}
\def\Id{{\rm Id}}
\def\spann{\mathop{\rm span}}
\def\fl{\mathop{\tt fl}}
\def\op{\mathop{\tt op}}
\def\cost{{\rm cost}}
\def\emac{\varepsilon_{\mathsf{mach}}}
\def\cond{\mathop{\mathsf{cond}}}
\def\card{{\rm card}}
\def\mod{{\rm mod}\;}
\def\trace{{\rm trace}}
\def\Prob{\mathop{\rm Prob}}
\def\length{{\rm length}}
\def\Diag{{\rm Diag}}
\def\diam{{\rm diam}}
\def\Tr{{\mathsf{Tr}}}
\def\arg{{\rm arg}\,}
\def\ot{\leftarrow}
\def\transp{^{\rm t}}
\def\bL{{\bf L}}
\def\bR{{\bf R}}
\def\bC{\mathbf{C}}
\def\bfE{{\bf E}}
\def\bd{{\bf d}}
\def\adj{{\rm adj}\:}
\def\cruz{\raise0.7pt\hbox{$\scriptstyle\times$}}
\def\todif{\stackrel{\scriptscriptstyle\not =}{\to}}
\def\nedif{\stackrel{\scriptscriptstyle\not=}{\nearrow}}
\def\sedif{\stackrel{\scriptscriptstyle\not=}{\searrow}}
\def\nodown{\downarrow\mkern-15.3mu{\raise1.3pt\hbox{$\scriptstyle\times$}}}
\def\noto{\to\mkern-20mu\cruz}
\def\notox{\to\mkern-23mu\cruz}
\def\alg{{\hbox{\scriptsize\rm alg}}}
\def\ttl{{\tt l}}
\def\tto{{\tt o}}
\def\ttml{\bar{\tt l}}
\def\ttmo{\bar{\tt o}}
\def\ri{R^\infty}
\def\Oh{{\cal O}}
\def\parts{{\cal P}}
\def\coeff{{\hbox{\rm coeff}}}
\def\Gr{{\hbox{\rm Gr}}}
\def\Error{{\hbox{\tt Error}}\,}
\def\appleq{\hbox{\lower3.5pt\hbox{$\;\:\stackrel{\textstyle<}{\sim}\;\:$}}}

\def\bolita{\scriptscriptstyle\bullet}


\def\JACM{Journal of the ACM}
\def\CACM{Communications of the ACM}
\def\ICALP{International Colloquium on Automata, Languages
            and Programming}
\def\STOC{annual ACM Symp. on the Theory
          of Computing}
\def\FOCS{annual IEEE Symp. on Foundations of Computer Science}
\def\SIAM{SIAM J. Comp.}
\def\SIOPT{SIAM J. Optim.}
\def\BSMF{Bulletin de la Soci\'et\'e Ma\-th\'e\-ma\-tique de France}
\def\CRAS{C. R. Acad. Sci. Paris}
\def\IPL{Information Processing Letters}
\def\TCS{Theoret. Comp. Sci.}
\def\BAMS{Bulletin of the Amer. Math. Soc.}
\def\TAMS{Transactions of the Amer. Math. Soc.}
\def\PAMS{Proceedings of the Amer. Math. Soc.}
\def\JAMS{Journal of the Amer. Math. Soc.}
\def\LNM{Lect. Notes in Math.}
\def\LNCS{Lect. Notes in Comp. Sci.}
\def\JSL{Journal for Symbolic Logic}
\def\JSC{Journal of Symbolic Computation}
\def\JCSS{J. Comput. System Sci.}
\def\JoC{J. Compl.}
\def\MP{Math. Program.}
\sloppy

\bibliographystyle{plain}


\def\CPRi{{\rm \#P}_{\kern-2pt\R}}

\newcommand{\bfb}{{\boldsymbol{b}}}
\newcommand{\bfd}{{\boldsymbol{d}}}
\newcommand{\bff}{{\boldsymbol{f}}}
\newcommand{\bfx}{{\boldsymbol{x}}}
\newcommand{\bfa}{{\boldsymbol{a}}}
\newcommand{\ol}[1]{\overline{#1}}
\def\msC{\mathop{\mathscr C}}
\def\mcC{{\mathcal C}}
\def\DD{{\mathcal D}}
\def\mcN{{\mathcal N}}
\def\mcG{{\mathcal G}}
\def\mcU{{\mathcal U}}
\def\mcM{{\mathcal M}}
\def\mcZ{{\mathcal Z}}
\def\mcB{{\mathcal B}}
\def\mcX{{\mathcal X}}
\def\PP{{\mathscr P}}
\def\AA{{\mathscr A}}
\def\scC{{\mathscr C}}
\def\scO{{\mathscr O}}
\def\scE{{\mathscr E}}
\def\sG{{\mathscr G}}
\def\mZ{{\mathcal Z}}
\def\mI{{\mathcal I}}
\def\sH{{\mathscr H}}
\def\sF{{\mathscr F}}
\def\bD{{\mathbf D}}
\def\nrr{{\#_{\R}}}
\def\Sd{{\Sigma_d}}
\def\oPp{\overline{P'_a}}
\def\oP{\overline{P_a}}
\def\oDH{\overline{DH_a^\dagger}}
\def\oH{\overline{H_a}}
\def\ua{\overline{\a}}
\def\Hd{\HH_{{\boldsymbol{d}}}}
\def\Hdm{\Hd[m]}
\def\Lg{{\rm Lg}}
\def\sfC{{\mathsf{C}}}
\def\Rk{{\operatorname{rank}}}

\def\bx{{\bf x}}
\def\ii{{\'{\i}}}

\def\P{\mathbb P}
\def\Exp{\mathop{\mathbb E}}

\newcommand{\binomial}[2]{\ensuremath{{\left(
\begin{array}{c} #1 \\ #2 \end{array} \right)}}}

\newcommand{\HH}{\ensuremath{\mathcal H}}
\newcommand{\diag}{\mathbf{diag}}
\newcommand{\CH}{\mathsf{CH}}
\newcommand{\Cone}{\mathsf{Cone}}
\newcommand{\SCH}{\mathsf{SCH}}


\newcommand{\NPR}{\mathsf{NP}_{\R}}
\newcommand{\NPc}{\mathsf{NP}}
\newcommand{\sfH}{\mathsf{H}}
\newcommand{\SPR}{\#\mathsf{P}_{\R}}
\newcommand{\PSPACE}{\mathsf{PSPACE}}
\newcounter{line}

\newcommand{\macheps}{\varepsilon_{\mathrm{m}}}
\newcommand{\sgn}{\mathrm{sgn}}

\begin{title}
{\LARGE {\bf Computing the Homology of Real Projective Sets}}
\end{title}
\author{Felipe Cucker\thanks{Partially supported by a GRF grant
from the Research Grants Council of the Hong Kong
SAR (project number CityU  11310716).}
\\
Dept. of Mathematics\\
City University of Hong Kong\\
HONG KONG\\
e-mail: {\tt macucker@cityu.edu.hk}
\and
Teresa Krick\thanks{Partially supported by grants
BID-PICT-2013-0294, UBACyT-2014-2017-20020130100143BA  and PIP-CONICET 2014-2016-11220130100073CO.}\\
Departamento de Matem\'atica \& IMAS\\
Univ. de Buenos Aires \&\ CONICET\\
ARGENTINA\\
e-mail: {\tt krick@dm.uba.ar}
\and
Michael Shub\\
Department of Mathematics\\ City College and the Graduate Center of CUNY\\
New York\\
USA\\
e-mail: {\tt  mshub@ccny.cuny.edu}
}

\date{}
\makeatletter
\maketitle
\makeatother

\begin{quote}
{\small
{\bf Abstract.}
We describe and analyze a numerical algorithm for computing the
homology (Betti numbers and torsion coefficients) of real projective
varieties. Here numerical means that the algorithm is numerically
stable (in a sense to be made precise). Its cost depends on the condition
of the input  as well as on its size  and is singly exponential in the
number of variables (the dimension of the ambient space) and
polynomial in the condition and the degrees of the defining
polynomials. In addition, we show that outside of an exceptional
set of measure exponentially small in the size of the data,
the algorithm takes exponential time.
{\bf AMS classification numbers:} 65Y20, 65H10, 55U10.
}\end{quote}

\section{Introduction}

This paper describes and analyzes, both in terms of complexity and
numerical stability, an algorithm to compute the topology of a
real projective set.

The geometry of the sets of zeros of polynomials equalities, or more
generally solutions of polynomial inequalities, is strongly tied to
complexity theory. The problem of deciding whether such a set is
nonempty is the paramount $\NPR$-complete problem
(i.e., $\NPc$-complete over the reals)~\cite{BSS89}; deciding
whether it is unbounded is
$\sfH\exists$-complete and whether a point is isolated on it is
$\sfH\forall$-complete~\cite{BC09}; computing its Euler
characteristic, or counting its points (in the zero dimensional case),
$\SPR$-complete~\cite{BC03}, \dots

We do not describe complexity classes in these pages. We content
ourselves with the observation that such classes are characterized by
restrictions in the use of specific resources (such as computing time
or working space)
and that complete problems are representatives for
them.  In this sense, the landscape of classes demanding an increasing
amount of resources is paralleled by a collection of problems whose
solution appears to be increasingly difficult.

Among the problems whose complexity is poorly understood, the
computation of the homology of algebraic or semialgebraic sets ---and
by this we mean the computation of all their Betti numbers and torsion
coefficients--- stands out. The use of Cylindrical Algebraic
Decomposition~\cite{Collins,Wut76}
allows one to compute a triangulation of the set at hand
(and from it, its homology) with a running time doubly exponential in
the number of variables (the dimension of the ambient
space). On the other hand, the $\SPR$-hardness of
computing the Euler characteristic (a simpler problem) mentioned above
or the $\PSPACE$-hardness of the problem of computing all Betti
numbers of a complex algebraic (or projective) set defined over $\Z$,
see~\cite{Scheib:07}, make clear that the existence of
subexponential algorithms for the computation of the
homology is unlikely. The obvious question is whether exponential
time algorithms for this task exist.

A number of results in recent years have made substantial progress
towards an answer to this question.
Saugata Basu and collaborators provide algorithms computing the first
Betti number of a semialgebraic set in single exponential time
(an algorithm to compute the zeroth Betti number within these bounds
was already known)~\cite{BaPoRo:08}, as well as an algorithm
computing the top $\ell$ Betti numbers with cost doubly exponential
in $\ell$ (but polynomial for fixed $\ell$)~\cite{Basu:08}. More recently,
Peter Scheiblechner~\cite{Scheib:12} considered
the class of smooth complex projective varieties and exhibited
an algorithm computing all the Betti numbers
(but not the torsion coefficients as the paper actually computes
the de Rham homology) for sets in this
class in single exponential time.

All the algorithms mentioned above are ``symbolic'', they are direct
(as opposed to iterative) and are not meant to work under finite
precision. Actually, numerical instability has been observed for many
of them and very recent results~\cite{NoTo:15} give some theoretical
account for this instability.  And partly motivated by this observed
instability, an interest in numerical algorithms has developed in
tandem with that on symbolic algorithms. An example of the former that
bears on this paper is the algorithm in~\cite{CS98} to decide
feasibility of semialgebraic sets. The idea was to decide the
existence of the desired solution by exploring a grid. While this grid
would have exponentially many points, the computation performed at
each such point would be fast and accurate, thus ensuring numerical
stability in the presence of round-off errors. Both the running time of the
algorithm (directly related to the size of the grid) and the machine
precision needed to ensure the output's
correctness, were shown to depend on a
condition number for the system of polynomial inequalities defining
the semialgebraic set at hand.

These ideas were extended in~\cite{CKMW1,CKMW2,CKMW3}
to describe and analyze a numerical algorithm for the more difficult
question of counting points in zero-dimensional projective sets.
Note that in this case
the number to be computed coincides with the zeroth
Betti number of the set (number of connected components),
while higher Betti numbers are all zero.

We now extend them once more to solve the (even more difficult)
problem of computing all the homology groups for
projective (or spherical) algebraic sets.
\medskip

In order to state our result, we need to introduce some notation.

Let $m\le n$, $d_1,\ldots,d_m\in \N$ and $\bfd = (d_1, \dots, d_m)$.
We will denote by $\Hdm$ the space of polynomial systems
$f=(f_1,\ldots,f_m)$ with $f_i\in\R[X_0,\ldots,X_n]$ homogeneous of
degree~$d_i$. We may assume here that $d_i\ge 2$ for $1\le i\le m$,
since otherwise we could reduce the input to a system with fewer
equations and unknowns. We set $D:=\max\{d_i,\,1\le i\le m\}$ and
$N:=\dim_{\R}\Hdm=\sum_{i=1}^m \binom{n+d_i}{n}$. Note that the
last is the {\em size} of the system $f$ in the sense that it is the
number of reals needed to specify this system.

We associate to
$f\in \Hdm$ its zero sets $\mcM_{\IS}:=Z_{\IS^n}(f)$ on
the unit sphere $\IS^n\subset \R^{n+1}$ and
$\mcM_{\P}:=Z_{\P^n}(f)$ on the projective
space $\proj^n(\R)$. The former is the intersection of
the cone of zeros $\mcZ:=Z_{\R^{n+1}}(f)$ of $f$ in $\R^{n+1}$ with
$\IS^n$ and the latter is the quotient of $\mcM_{\IS}$ by identifying
antipodal points.
For a generic system $f$, both $\mcM_{\IS}$ and
$\mcM_{\P}$ are smooth manifolds of dimension $n-m$.
We also associate to $f$ a condition number
$\kappa(f)$ (whose precise definition will be
given in \S\ref{sec:cond-numb} below). Finally, we endow the linear
space $\Hd[m]$ with the Weyl inner product (also
defined in~\S\ref{sec:cond-numb}) and consider the unit
sphere $\IS^{N-1}\subset \Hd[m]$ with respect to the norm induced
by it.

\begin{theorem}\label{thm:Main}
We describe an algorithm that, given $f\in\Hd[m]$, returns the
Betti numbers and torsion coefficients of $\mcM_{\IS}$ (or
of $\mcM_{\P}$), with the following properties.
\begin{description}
\item[(i)]
Its cost $\cost(f)$ on input $f$ is bounded by
$(nD\kappa(f))^{\Oh(n^2)}$.
\item[(ii)]
Assume $\IS^{N-1}$ is endowed with the uniform probability
measure. Then, with probability at least $1-(nD)^{-n}$ we have
 $\cost(f)\leq (nD)^{\Oh(n^3)}$.
\item[(iii)]
Similarly, with probability at least $1-2^{-N}$ we have
 $\cost(f)\leq 2^{\Oh(N^2)}$.
\item[(iv)]
The algorithm is numerically stable.
\end{description}
\end{theorem}

We give the proof of Theorem~\ref{thm:Main} in several steps.
Part~(i) is shown in Propositions~\ref{prop:betti-S}
and~\ref{prop:betti-P}. Parts~(ii) and~(iii) are in
Corollary~\ref{corol:probs}. We devote Section~\ref{sec:stability}
to both define what we mean by numerical stability (in a context
where we are computing integer numbers) and to sketch why
our algorithm is numerically stable.

\begin{remark}
Parts~(ii) and~(iii) in the statement fit well within the setting of
{\em weak complexity analysis} recently proposed in~\cite{wacco}
(but see also~\cite[Theorem~4.4]{Kostlan88} for a
predecessor of this setting).
The idea here is to exclude from the analysis
a set of outliers of exponentially small measure (a probability
measure in the space of data is assumed). This exclusion may
lead to dramatic differences in the quantity to be bounded
and provide a better agreement between theoretical analysis
and computational experience. A case at hand, studied
in~\cite{wacco}, is that of the power method to compute dominant
eigenpairs. It is an algorithm experienced as efficient in practice
(say for symmetric or Hermitian matrices) but whose expected
number of iterations (for matrices drawn from the Gaussian orthogonal
or unitary ensembles, respectively) is known to be
infinite~\cite{Kostlan88}. Theorem~1.4 in~\cite{wacco} shows that
the expected number of iterations conditioned to excluding a set
of exponentially small measure is polynomially bounded in the
dimension $n$ of the input matrix. The authors call this form
of analysis {\em weak average-case}. Parts~(ii) and~(iii) in the statement
can be seen as a form of weak worst-case analysis establishing
weak worst-case exponential complexity.
\end{remark}

Our algorithm relies on an extension of the ideas in~\cite{CS98}
---the use of grids, an exclusion test, and the use of the
$\alpha$-theory of Smale to detect zeros of a polynomial system in
the vicinity of a point at hand--- to construct a covering of
$\mcM_{\IS}$ by open balls in $\R^{n+1}$ of the same radii. This common radius
is chosen to ensure that the union of the balls in the covering is
homotopically equivalent to $\mcM_{\IS}$. The Nerve Theorem
then ensures that this union is homotopically equivalent to the
nerve of the covering and we can compute the homology groups
of $\mcM_{\IS}$ by computing those of the said nerve. We explain
the basic ingredients (condition numbers, Smale's $\alpha$-theory,
the exclusion lemma, \dots) in Section~\ref{sec:basics}. Then, in
Section~\ref{sec:algorithm}, we describe and analyze the computation
of  the covering. Section~\ref{sec:betti} uses this covering to actually
compute the homology groups (part~(i) in Theorem~\ref{thm:Main})
and Section~\ref{sec:random} establishes the probability estimates
(parts~(ii) and~(iii) in Theorem~\ref{thm:Main}).
Section~\ref{sec:proofs} is devoted to prove a number of results
which, to allow for a streamlined exposition, were only stated in
Section~\ref{sec:basics}. One of them,
Theorem~\ref{prop:tau-gamma}, links the $\gamma$-invariant
of Smale with the injectivity radius $\tau(f)$ of the normal bundle
of $\mcM_{\IS}$ (in turn related to a number of metric properties of
algebraic spherical (or projective) sets). This connection is,
to the best of our knowledge, new and is interesting per se.
Finally, and as already mentioned, Section~\ref{sec:stability}
deals with issues of finite-precision and numerical stability.
\bigskip

\noindent
{\bf Acknowledgments.}
We are grateful to Peter B\"urgisser who
suggested the topic of this paper to us  and to the Simons Institute for receiving us in the Fall of 2014, which was where and when the suggestion was made. We also owe an anonymous referee for
his very precise and enlightening comments.

\section{The basic ingredients}\label{sec:basics}

\subsection{Condition numbers}\label{sec:cond-numb}

We need a condition number as a complexity (and accuracy)
parameter. To define one we first fix a norm on the space
$\Hdm$. We follow the (by now well-established) tradition of using
the Weyl norm, which is invariant under
the action of orthogonal transformations in $\R^{n+1}$:
for $f=(f_1,\ldots,f_m)$ with
$f_i=\sum_{|\bfa|=d}f_{i,\bfa} X^\bfa$, this is
$\|f_i\|^2= \sum _{|\bfa|=d}f_{i,\bfa}^2{d\choose \bfa}^{-1}$
and then $\|f\|^2 := \sum_{1\le i\le m}\|f_i\|^2$.
See e.g.~\cite[\S16.1]{Condition} for details.

For a point $\xi\in\R^{n+1}$ we denote by
$Df(\xi)=\Big( \dfrac{\partial f_i}{\partial x_j}(\xi)
 \Big)_{1\le i\le m, 0\le j\le n}:\R^{n+1}\to\R^m$
the derivative of $f$ at $\xi$. We also write
$$
  \Delta(\xi):={\small \left[ \begin{array}{ccc} \|\xi\|^{d_1-1}\sqrt{d_1} \\
  & \ddots &\\
  & &  \|\xi\|^{d_m-1}\sqrt{d_m} \end{array} \right]}
$$
(or simply $\Delta$, if $\xi\in\IS^n$).

The condition of $f$ at a zero $\xi\in\R^{n+1}\setminus\{0\}$ has been
well-studied in the series of
papers~\cite{Bez1,Bez2,Bez3,Bez4,Bez5}. It is defined as $\infty$ when
the derivative $Df(\xi)$ of $f$ at $\xi$ is not surjective, and when
$Df(\xi)$ is surjective as
\begin{equation}\label{eq:munorm}
   \mun(f,\xi):=\|f\| \big\|Df(\xi)^\dagger \Delta(\xi)\big\|,
\end{equation}
where
$Df(\xi)^\dagger: \R^{m} \to \R^{n+1}$ is the  {\em Moore-Penrose inverse}
of the full-rank matrix $Df(\xi)$, i.e.
$Df(\xi)^\dagger=Df(\xi)\transp(Df(\xi) \, Df(\xi)\transp)^{-1}$,
where $Df(\xi)\transp$ is the transpose of $Df(\xi)$. This coincides
 with the inverse of the restricted linear map
$Df(\xi)|_{(\ker Df(\xi))^\perp}$.
Also, the norm in $\|Df(\xi)^\dagger \Delta(\xi)\|$ is the spectral norm.

Since  the expression in the right of~\eqref{eq:munorm}
is well-defined for arbitrary points $x\in\IS^n$, we can
define $\mun(f,x)$ for any such point.

For 0-dimensional homogeneous systems, that is, for systems
$f\in\Hd[n]$, the quantity $\mun(f,x)$ in~\eqref{eq:munorm}
is occasionally defined
differently,  by replacing $Df(x)^{\dagger}$
by $(Df(x)_{|T_x})^{-1}$. Here $T_x$
denotes the orthogonal complement of $x$ in $\R^{n+1}$
and we are inverting the restriction of the derivative $Df(x)$
to this space (see~\cite[\S16.7]{Condition}).
This definition only makes sense when $m=n$ as in this case
the restriction $(Df(x)_{|T_x})^{-1}:T_x\to\R^n$ is a linear map
between spaces of the same dimension. This is not the case
when $m\neq n$. Hence the
use here of the Moore-Penrose derivative.

To define the condition of a system $f$ it is not enough to just
consider the condition at its zeros. For points $x\in\R^{n+1}$
where $\|f(x)\|$ is non-zero but small, small perturbations of $f$
can turn $x$ into a new zero (and thus change the topology
of $\mcZ$). Following an idea going back
to~\cite{Cucker99b} and developed in this context
in~\cite{CKMW3} we define
$$
\kappa(f,x):=
    \frac{\|f\|}{\big\{\|f\|^2\mun^{-2}(f,x)
     +\|f(x)\|^2\big\}^{1/2}}
$$
where $\mun(f,x)$ is defined as in~\eqref{eq:munorm} for $x\in \IS^n$,
with the convention that $\infty^{-1}=0$ and $0^{-1}=\infty$, and
\begin{equation}\label{eq:kappa}
\kappa(f):=\max_{x\in\IS^n}\kappa(f,x).
\end{equation}

\begin{remark}\label{rem:inv}
For any $\lambda\neq 0$ we have
$\mun(f,x)=\mun(f,\lambda x)$, since when
$Df(x)$ is surjective, $Df(\lambda x)^\dagger= \Big(\Lambda
Df(x)\Big)^\dagger = Df(x)^\dagger \Lambda^{-1} $ for
$\Lambda={\scriptsize \left[\begin{array}{ccc}\lambda^{d_1-1} & & \\ &
\ddots & \\ & & \lambda^{d_m-1}\end{array} \right]} $. Similarly,
$\mun(f,\xi)=\mun(\lambda f, \xi)$ for all $\lambda\neq 0$, and
consequently, $\kappa(\lambda f)=\kappa(f)$.
\end{remark}

Note that $\kappa(f)=\infty$ if only if there exists $\xi \in \IS^n$
such that $f(\xi)=0$ (i.e $\xi\in \mcM_\IS$) and $Df(\xi)$ is not
surjective, i.e., $f$ belongs to the {\em set of ill-posed systems}
\begin{equation}\label{def:sigma}
      \Sigma_\R:=\big\{f\in\Hdm\mid \mbox{$\exists \, \xi\in\IS^n$ such that
    $f(\xi)=0$ and $\Rk( Df(\xi))<m$}\big\}.
\end{equation}

The following result is proved in Section~\ref{sec:kappa-dist}.
It extends a statement originally shown for square systems
in~\cite{CKMW2} (see also~\cite[Theorem~19.3]{Condition}).

\begin{proposition}\label{prop:kappa-dist}
For all $f\in\Hdm$,
$$
  \frac{\|f\|}{{\sqrt{2}}\,\dist(f,\Sigma_\R)} \le   \kappa(f)
  \le \frac{\|f\|}{\dist(f,\Sigma_\R)}.
$$
\end{proposition}

We prove the following in Section~\ref{sec:lipschitz2}.

\begin{proposition}\label{prop:lipschitz2}
Let $m\leq n+1$.
For all $f\in\Hd[m]$, $0\le \e \le \frac12$  and $y,z\in\IS^n$
such that
$$
     \|y-z\|\leq\frac{2\e}{D^{3/2}\mun(f,y)}
$$
we have
$$
   \frac1{1+\frac52\e}\,\mun(f,y) \leq
   \mun(f,z) \leq \Big(1+\frac52\e\Big)\mun(f,y).
$$
\end{proposition}

\subsection{Moore-Penrose Newton and point estimates}

Let $f:\R^{n+1}\to\R^m$, $m\leq n+1$, be analytic.
The {\em Moore-Penrose Newton operator} of
$f$ at $x\in \R^{n+1}$ is defined
(see~\cite{AllGeo90}) as
$$
    N_f(x):= x-Df(x)^\dagger f(x).
$$
We say that it is well-defined if $Df(x)$ is surjective.

\begin{definition}
Let $x \in \R^{n+1}$. We say that  $x$ {\em converges to a zero of $f$} if the sequence $(x_k)_{k\ge 0}$  defined as $x_0:= x$ and
$x_{k+1}:=N_f(x_k)$ for $k\ge 0$ is well-defined and converges to a zero of $f$.
\end{definition}

Following ideas introduced by Steve Smale in~\cite{Smale86},
the following three quantities were associated to a point
$x\in \R^{n+1}$ in~\cite{Bez4},
\begin{eqnarray*}
\b(f,x)&:=&\|Df(x)^\dagger \,f(x)\||\\
\gamma(f,x)&:=&\max_{k>1}
\left\|Df(x)^\dagger \frac{D^kf(x)}{k!}\right\|^{\frac{1}{k-1}}\\
\a(f,x)&:=&\b(f,x)\gamma(f,x),
\end{eqnarray*}
when $Df(x)$ is surjective,
and $\alpha(f,x)=\beta(f,x)=\gamma(f,x)=\infty$ when
$Df(x)$ is not surjective.
The quantity $\b(f,x)=\|N_f(x)-x\|$ measures the length of the Newton
step at $x$. The value of $\gamma(f,\xi)$, at a zero $\xi$ of $f$,
is related to the radius of the neighborhood of points that converge to the zero $\xi$
of $f$, and the meaning of $\a(f,x)$
is made clear in the main theorem in the theory of point
estimates.

\begin{theorem}\label{alpha-th}
Let $f:\R^{n+1}\to\R^m$, $m\leq n+1$, be analytic.
Set $\a_0=0.125$. Let $x \in \R^{n+1}$ with
 $\a(f, x)< \a_0$, then $x$ converges to a zero $\xi$ of $f$ and
$\|x-\xi \|<2\beta(f,x)$.
Furthermore, if $n+1=m$ and $\a(f,x)\leq 0.03$, then all
points in the ball of center $x$ and radius
$\frac{0.05}{\gamma(f,x)}$ converge to the same zero of $f$.
\end{theorem}

\proof
In~\cite[Th.~1.4]{Bez4} it is shown that under the stated hypothesis,   $x$  converges to a zero $\xi$ of $f$ and
$$
    \|x_{k+1}-x_k\|\le
    \left(\frac{1}{2}\right)^{2^k-1}\|x_1-x_0\|= \left(\frac{1}{2}\right)^{2^k-1}\beta(f,x).
$$
Therefore $$\|x_{i+1}-x \|\le \sum_{0\le k\le i}\left(\frac{1}{2}\right)^{2^k-1}\beta(f,x) < (2-\frac{1}{8})\beta(f,x).$$ This implies the first statement.
The second is Theorem~4 and Remarks~5, 6 and 7
in~\cite[Ch.~8]{BCSS98}.
\eproof

In what follows we will apply the theory of point estimates
to the case of polynomial maps $f=(f_1,\ldots,f_m)$.
In the particular case where the $f_i$ are homogeneous,
the invariants $\a,\b$ and $\gamma$ are themselves
homogeneous in $x$. We have
$\beta(f,\lambda x)=\lambda\beta(f,x)$,
$\gamma(f,\lambda x)=\lambda^{-1}\gamma(f,x)$, and
$\a(f,\lambda x)=\a(f,x)$, for all $\lambda\neq 0$. This
property motivates the following projective version for
them:
\begin{eqnarray*}
\b_{\hmg}(f,x)&:=&\|x\|^{-1}\|Df(x)^\dagger \,f(x)\|\\
\gamma_{\hmg}(f,x)&:=&\|x\|\max_{k>1}
\left\|Df(x)^\dagger \frac{D^kf(x)}{k!}\right\|^{\frac{1}{k-1}}\\
\a_{\hmg}(f,x)&:=&\b_{\hmg}(f,x)\gamma_{\hmg}(f,x),
\end{eqnarray*}
These projective versions coincide with the previous expressions
when $x\in \IS^n$ and an $\a$-Theorem for them is easily derived
from Theorem~\ref{alpha-th}  above. Furthermore,
$\b_{\hmg}$ still measures the (scaled) length of the Newton step,
and $\gamma_{\hmg}$ relates to the condition number
via the following bound
(known as the Higher Derivative Estimate),
\begin{equation}\label{eq:HDE}
\gamma_{\hmg}(f,x)\leq \frac12D^{3/2}\mun(f,x).
\end{equation}
The proof is exactly the one of~\cite[Th.~2, p.~267]{BCSS98} which
still holds for $m\le n$ and $Df(x)^\dagger$ instead of
$Df(x)|_{T_x}^{-1}$.

We now  move to  ``easily computable''
versions $\ba,\bb$ and $\bg$, which we define
for $x\in\IS^n$:
\begin{eqnarray}\label{eq:abc}
\bb(f,x)&:=&\mun(f,x)\frac{\|f(x)\|}{\|f\|} \nonumber\\
\bg(f,x)&:=&\frac12D^{3/2}\mun(f,x)\\
\ba(f,x)&:=&\bb(f,x)\bg(f,x)=\frac12 D^{3/2}\mun^2(f,x)\frac{\|f(x)\|}{\|f\|}.\nonumber
\end{eqnarray}

For $x\in \IS^n$, \eqref{eq:HDE} therefore says that $\gamma(f,x)\le \bg(f,x)$. We also observe that $\beta(f,x)\le \bar\beta(f,x)$ since
$$
  \beta(f,x)=
  \left\|Df(x)^\dagger f(x)\right\|\le
  \left\|Df(x)^\dagger\right\|\|f(x)\|\le
  \|f\|\|Df(x)^\dagger \Delta\|\frac{\|f(x)\|}{\|f\|}=\bar\beta(f,x).
$$
Therefore $\alpha(f,x)\le \bar \alpha(f,x)$.

\subsection{Curvature and coverings}

A crucial ingredient in our development is a result
in a paper by Niyogi, Smale  and
Weinberger~\cite[Prop.7.1]{NiSmWe08}. The context
of that paper (learning on manifolds) is different from
ours but this particular result, linking curvature and
coverings, is, as we said, central to us.

Consider a compact Riemannian submanifold $\mcM$ of
a Euclidean space $\R^{n+1}$. Consider as well a finite
collection of points $\mcX=\{x_1,\ldots,x_K\}$ in $\R^{n+1}$
and also $\e>0$. We are interested in conditions guaranteeing
that the union of the open balls
$$
    U_\e(\mcX):= \bigcup_{x\in\mcX} B(x,\e)
$$
covers $\mcM$ and is homotopically equivalent to it. These
conditions involve two notions which we next define.

We denote by $\tau(\mcM)$ the {\em injectivity radius of the
normal bundle of $\mcM$}, i.e., the largest $t$
such that the open normal bundle around $\mcM$ of radius $t$
$$
   N_t(\mcM):=\big\{(x,v)\in\mcM\times\R^{n+1}\mid v\in N_x\mcM,
               \|v\|<t\big\}
$$
is embedded in $\R^{n+1}$. That is, the largest $t$ for which
$\phi_t:N_t(\mcM)\to\R^{n+1}$, $(x,v)\mapsto x+v$, is injective.
Therefore,
its image $\Tub_{\tau(\mcM)}$ is an open tubular neighborhood
of $\mcM$ with its canonical orthogonal projection map
$\pi_0:\Tub_{\tau(\mcM)}\to \mcM$ mapping every point
$x\in\Tub_{\tau(\mcM)}$ to the (unique) point in $\mcM$
closest to $x$.  In particular, $\mcM$ is
a deformation retract of $\Tub_{\tau(\mcM)}$.

Also, we recall that the Hausdorff distance between two
subsets $A,B\subset\R^{n+1}$ is defined as
$$
   d_H(A,B):=\max \Big\{\sup_{a\in A}\inf_{b\in B} \|a-b\|,
   \sup_{b\in B}\inf_{a\in A} \|a-b\|\Big\}.
$$
If both $A$ and $B$ are compact, we have that $d_H(A,B)\le r$
if and only if for all $a\in A$ there exists $b\in B$ such that $\|a-b\|\leq r$
and for all $b\in B$ there exists $a\in A$ such that $\|a-b\|\leq r$.

The following is a slight variation of~\cite[Prop.7.1]{NiSmWe08}.

\begin{proposition}\label{prop:SNW}
Let $\bt \le  \tau(\mcM)$ and $0<r<(3-\sqrt{8})\bt$.
If $d_H(\mcX,\mcM)\le r$
then $\mcM$ is a deformation retract of $U_\e(\mcX)$
for every $\e$ satisfying
\begin{equation}\tag*{\qed}
  \e\in\left(\frac{(r+\bt)-\sqrt{r^2+\bt^2-6r\bt}}{2},
    \frac{(r+\bt)+\sqrt{r^2+\bt^2-6r\bt}}{2}\right).
\end{equation}
\end{proposition}

\begin{remark}\label{rem:tau}
If we start with  $r>0$ for which $6r<\tau(\mcM)$ we can take
$\bt:=6r$. In this case the interval we obtain for the admissible values of
$\e$ is $[3r,4r]$.
\end{remark}

The quantity $\tau(\mcM)$ is strongly related to the curvature of
$\mcM$ as shown in Propositions~6.1, 6.2, and 6.3
in~\cite{NiSmWe08}. Even though we won't make use of these
results, we summarize them in the following statement.

\begin{theorem}\label{th:tau}
Let $\tau:=\tau(\mcM)$.

{\bf (i)\quad} The norm of the second fundamental form of $\mcM$
is bounded by $\frac1{\tau}$ in all directions.

{\bf (ii)\quad}For $p,q\in\mcM$ let $\phi(p,q)$ be the angle between
their tangent spaces $T_p$ and $T_q$, and $d_{\mcM}(p,q)$ their
geodesic distance. Then
$\cos(\phi(p,q))\geq 1-\frac{1}{\tau}d_{\mcM}(p,q)$.

{\bf (iii)\quad}For $p,q\in\mcM$,
$\displaystyle d_{\mcM}(p,q) \leq \tau -\tau
\sqrt{1-\frac{2\|p-q\|}{\tau}}$. \eproof
\end{theorem}

\subsection{Curvature and condition}

Theorem~\ref{th:tau} shows a deep relationship
between the curvature of a submanifold $\mcM$ of
Euclidean space and the value of $\tau(\mcM)$. One
of the main results in this paper is a further connnection,
for the particular case where $\mcM=\mcM_{\IS}$,
the set of zeros of $f\in\Hdm$ in $\IS^{n}$, between
$\tau(\mcM_{\IS})$ and the values of $\gamma$ on $\mcM_{\IS}$.
Define
$$
 \tau(f):=\tau(  \mcM_\IS) \quad \mbox{and} \quad
 \Gamma(f):=\max_{x\in\mcM_{\IS}} \max\{1,\gamma(f,x)\}.
$$
In Section~\ref{sec:tau-gamma} we prove the following.

\begin{theorem}\label{prop:tau-gamma}
We have
$$
   \tau(f)\geq \frac{1}{87\,\Gamma(f)}.
$$
\end{theorem}

Note that as $\max\{1,\gamma(f,x)\}\leq \bg(f,x)$ we
obtain
\begin{corollary}\label{gammasomb}$$
   \tau(f)\geq \frac{1}{87\,\overline\Gamma(f)}.
$$\end{corollary}
where $\overline{\Gamma}(f):=\max_{x\in\mcM_{\IS}}\bg(f,x)$.

\subsection{Grids and exclusion results}
\label{sec:grids}

Our algorithm works on a grid $\mcG_\eta$ on $\IS^n$,
which we construct by
projecting onto $\IS^n$ a grid on the cube
$\sfC^n=\{y\in \R^{n+1}\mid \|y\|_\infty=1\}$.
We make use of the (easy to compute) bijections
$\phi:\sfC^n\to \IS^n$ and $\phi^{-1}:\IS^n\to \sfC^n$ given by
$\phi(y)=\frac{y}{\|y\|}$ and $\phi^{-1}(x)=\frac{x}{\|x\|_\infty}$.

Given $\eta:=2^{-k}$ for some $k\geq 1$, we consider the uniform
grid $\mcU_\eta$ of mesh $\eta$ on $\sfC^n$. This is the set of
points in $\sfC^n$ whose coordinates are of the form $i2^{-k}$ for
$i\in\{-2^k, -2^k+1,\ldots,2^k\}$, with at least one coordinate
equal to 1 or $-1$. We denote by $\mcG_\eta$ its image
by $\phi$ in $\IS^n$. An argument in elementary geometry shows that
for $y_1,y_2\in \sfC^n$,
\begin{equation}\label{eq:3dists}
  \|\phi(y_1)-\phi(y_2)\|\leq d_{\IS}(\phi(y_1),\phi(y_2))
  \leq \frac{\pi}{2}\|y_1-y_2\|
  \leq\frac{\pi}{2}\sqrt{n+1}\,\|y_1-y_2\|_\infty,
\end{equation}
where
$d_{\IS}(x_1,x_2):=\mbox{arccos}(\langle x_1,x_2\rangle)\in [0,\pi]$
denotes the angular distance, for $x_1,x_2\in\IS^n$.\\

Given $\varepsilon >0$, we denote by
$B(x,\varepsilon):=\{y\in \R^{n+1}\,|\, \|y-x\|<\varepsilon\}$,
for $x\in \R^{n+1}$, the open ball with respect to the Euclidean
distance, and by
$B_{\IS}(x,\varepsilon)=\{y\in \IS^n\,|\, d_\IS(y,x)<\varepsilon\}$,
for $x\in \IS^n$, the open ball with respect to the angular distance.
We also set from now on
 \begin{equation}\label{eq:sep}
\sep(\eta):=\eta\,\sqrt{n+1}\quad \mbox{and}
\quad \delta(f,\eta):=1.1 \sqrt{D(n+1)}\|f\|\eta.
 \end{equation}

\begin{lemma} \label{lem:cover}
The union
$\cup_{x\in\mcG_\eta} B(x,\sep(\eta))$ covers the sphere $\IS^n$.
\end{lemma}

\proof
Let $z\in\IS^n$ and $y=\phi^{-1}(z)\in\sfC^n$. There exists
$y'\in\mcU_{\eta}$ such that $\|y'-y\|_\infty\leq\frac{\eta}{2}$.
Let $x=\phi(y')\in\mcG_{\eta}$. Then, equation~\eqref{eq:3dists}
shows that $\|x-z\|\leq \frac{\eta}{2}\,\frac{\pi}{2}\,\sqrt{n+1}
< \eta\,\sqrt{n+1}$.
\eproof

In    \cite[Lem.~3.1]{CKMW1} and \cite[Lem.~19.22]{Condition},
the following Exclusion Lemma
is proved (the statement there is for $n=m$ but the proof holds
for general $m$).

\begin{lemma}{\bf (Exclusion lemma.)}\label{exclusionlemma}
Let $f \in \Hd[m]$ and $x, y \in \IS^n$ be such that
$0<d_{\IS}(x,y)\le \sqrt 2$. Then,
\[
  \|f(x)-f(y)\| < \|f\| \sqrt{D} \ d_{\IS}(x,y) .
\]
In particular, if $f(x) \ne 0$, there is no zero of $f$ in
the ball $B_{\IS}\big(x,\frac{\|f(x)\|}{\|f\|\sqrt{D}}\big)$. \eproof
\end{lemma}

\begin{corollary}\label{excl}
Let $\eta$ be such that $\sep(\eta)\leq\frac12$, and let
$x\in\IS^n$ satisfy $\|f(x)\|>\delta(f,\eta)$.
Then $f(y)\neq 0$ on the ball $B(x,\sep(\eta))$.
\end{corollary}

\proof
Let $y\in\R^{n+1}$ such that $\|y-x\|<\sep(\eta)\le \frac12$.
Define $h(\varepsilon)=\sqrt{2- 2\sqrt{1-\varepsilon^2}}$.
We have  $\|\phi(y)-x\|\le h(\|y-x\|)$.
Since $h(\varepsilon)/\varepsilon$ is monotonically increasing
on $[0,1]$,
$$
     \|\phi(y)-x\|\le 2 h(1/2) \|y-x\|<1.035 \|y-x\|< 0.5175
     \mbox{ for  } \ \|y-x\|<\frac12.
$$
Then,
$$
  d_{\IS}(\phi(y),x)= 2\arcsin\Big(\frac{\|\phi(y)-x\|}{2}\Big)
  \leq 1.012 \|\phi(y)-x\| < 1.1 \|x-y\|<1.1\, \sep(\eta)
$$
since $\arcsin$ is a convex function on the interval
$[0,0.5175]$.
Therefore the hypothesis on $\|f(x)\|$ implies that
$$
   \|f(x)\|>1.1\, \|f\|\sqrt{D}\,\sep(\eta) > \|f\|\sqrt{D}\,d_{\IS}(\phi(y),x)
$$
i.e.,  that $d_{\IS}(\phi(y),x)< \frac{\|f(x)\|}{\|f\|\sqrt{D}}$.
Lemma~\ref{exclusionlemma} then
shows, since $f(x)\ne 0$, that $f(\phi(y))\neq 0$  and we conclude that
$f(y)\neq0$ as $f$ is homogeneous.
\eproof

\section{Computing a homotopically equivalent covering}
\label{sec:algorithm}

Set $k:=\lceil\log_2 4\sqrt{n+1}\rceil$ so that $\sep(\eta)\leq
\frac14$ for $\eta=2^{-k}$, where $\sep(\eta)$ is defined
in~\eqref{eq:sep}. Our algorithm works on the grid $\mcG_\eta$
on $\IS^n$
constructed in the previous section, and makes use of the quantities
$\bb,\bg$ and $\ba$ introduced in~\eqref{eq:abc} and $\delta(f, \eta)$
defined in~\eqref{eq:sep}. We recall $\alpha_0:=0.125$.

\bigskip

\algoritmo
\begin{algorithm}\label{alg:covering}
{\sf Covering}\\
\inputalg{$f\in\Hdm$}\\
\specalg{$f\neq0$}\\
\bodyalg{
let $\eta:=2^{-k}$\\
repeat \\
\espacio $\mcX:=\emptyset$\\
\espacio $r:=\sqrt{\sep(\eta)}$\\
\espacio $\e:=3.5\,r$\\
\espacio
for all $x\in\mcG_\eta$ \\
\eespacio if $\ba(f,x)\leq\a_0$ and $\frac{1}{531\,\bg(f,x)}\geq r$
   and $2.2\,\bb(f,x)< r$ then\\
\eeespacio $\mcX:=\mcX\cup \{x\}$\\
\eespacio elsif $\|f(x)\|\geq \delta(f,\eta)$ then do nothing\\
\eespacio elsif go to (*)\\
\espacio return the pair $\{\mcX,\e\}$ and halt\\
\espacio end for\\
\espacio (*) $\eta:=\eta/2$\\
}
\Output{$\{\mcX,\e\}$}\\
\postcond{The algorithm halts if $f\not\in\Sigma_{\R}$. If
$\mcX=\emptyset$ then $\mcM_{\IS}$ is empty. Otherwise,
the set $\mcX$ is closed by the involution $x\mapsto -x$, and
the union of the balls $\{B(x,\e)\mid x\in\mcX\}$ covers
$\mcM_{\IS}$ and is homotopically equivalent to it.}
\end{algorithm}
\falgoritmo

In the sequel we  use the
quantity
\begin{equation}\label{eq:C_4}
   \bC:=\max\bigg\{12\,(n+1)D,\,
         \frac{531^2}{2}\sqrt{n+1} \,D^{3}\bigg\}.
\end{equation}
Note that we have $\bC=\Oh(n\,D^3)$.

\begin{proposition}\label{prop:algo-covering}
Algorithm {\sf Covering} is correct (it computes a list
$\{\mcX,\e\}$ satisfying its postconditions). Furthermore,
its cost is bounded by
$$
   \Oh\left(\log_2 (\bC\kappa(f))\,nN(2\bC\kappa^2(f))^n\right)
   =\left(nD\kappa(f)\right)^{\Oh(n)}
$$
and the number $K$ of points in the returned $\mcX$ is bounded by
$\left(nD\kappa(f)\right)^{\Oh(n)}$.
\end{proposition}

The rest of this section is devoted to prove
Proposition~\ref{prop:algo-covering}.

\begin{lemma}\label{lem:zeros}
Let $x\in\IS^n$ and $y\in\mcZ(f)$ be such that $\|x-y\|\leq 0.7$.
Then the point $\phi(y):=\frac{y}{\|y\|}\in\mcM_{\IS}$ satisfies
$\|x-\phi(y)\|\leq 1.1\|x-y\|$.
\end{lemma}

\proof The proof goes exactly as the proof of Corollary~\ref{excl}.
\eproof

The following two lemmas deal with the correctness of the algorithm.

Assume the algorithm halts for a certain value $\eta$.  Let
$\mcX$ be the set constructed  by the execution at this stage
and set $r=\sqrt{\sep(\eta)}$.

\begin{lemma}\label{lem:density}
The sets $\mcX$ and $\mcM_\IS$ satisfy $d_H(\mcX,\mcM_\IS)\le r$.
Furthermore, for all $y\in \mcM_{\IS}$,
there exists $x\in \mcX$ such that $\|y-x\|\le r^2$.
\end{lemma}

\proof
The points in $\mcG_\eta$ divide into two groups that
satisfy, respectively:
\smallskip

\noindent\fbox{$x\in\mcG_\eta\setminus \mcX$} This happens
when $\|f(x)\|\geq \delta(f,\eta)$,
and therefore, by Corollary~\ref{excl}, there are no
zeros of $f$ in the ball $B(x,\sep(\eta))=B(x,r^2)$.
\smallskip

\noindent\fbox{$x\in\mcX$} This happens when in particular
$\bar\alpha(f,x)<\alpha_0$,
and therefore,  by
Theorem~\ref{alpha-th}, there exist zeros of $f$ in the ball
$B(x,2\beta(f,x))\subset B(x,r/1.1)$ since $2.2\bar\beta(f,x)<r$.
This implies, because of Lemma~\ref{lem:zeros},
that $\mcM_{\IS}\cap B(x,r)\neq \emptyset$.
\smallskip

This last sentence shows that for $x\in \mcX$, there exists
$y\in \mcM_{\IS}$ with $\|y-x\|<r$.
In addition, since by Lemma~\ref{lem:cover},
$\cup_{x\in\mcG_\eta} B(x,r^2)$ covers the sphere
$\IS^n$ and there are no points of $\mcM_{\IS}$ in
$\cup_{x\in\mcG_\eta\setminus\mcX} B(x,r^2)$, it follows that
$\mcM_{\IS}\subset\cup_{x\in\mcX} B(x,r^2)$ and therefore for all $y\in \mcM_{\IS}$,
there exists $x\in \mcX$ such that $\|y-x\|\le r^2<r$.
This shows that $d_H(\mcX,\mcM_{\IS})\le r$.
\eproof

\begin{lemma}\label{lem:tau-r}
Let $\bt:=6r$. Then $\bt<\tau(f)$.
\end{lemma}

\proof
Let $y\in\mcM_{\IS}$ be such that $\overline{\Gamma}(f)=\bg(f,y)$,
for $\overline{\Gamma}(f)$
defined in Identity~\eqref{gammasomb}.
By Lemma~\ref{lem:density}
there exists $x\in\mcX$ such that $\|x-y\|< r$. Hence,
\begin{equation*}
 \|x-y\| \,<\, r\,\leq\,
 \frac{1}{531\,\bg(f,x)}
 \,=\,\frac{2}{531\,D^{3/2}\mun(f,x)}.
\end{equation*}
By Proposition~\ref{prop:lipschitz2}
(with $\e=\frac1{531}$) we  have
$\mun(f,y)\leq (1+\frac{5}{1062})\mun(f,x)
\leq 1.005\,\mun(f,x)$.
Consequently, $\bg(f,y)\le 1.005\bg(f,x)$ and therefore,
$$
  \bt = 6r \leq \frac{6}{531\,\bg(f,x)} \leq
       \frac{6.03}{531\,\bg(f,y)} < \frac{1}{87\,\bg(f,y)}
   = \frac{1}{87\,\overline{\Gamma}(f)} \leq \tau(f),
$$
the last by Theorem~\ref{prop:tau-gamma}.
\eproof

To bound the complexity we  rely on the following.

\begin{lemma}\label{lem:L2} Let $\bC$ be defined in \eqref{eq:C_4}.
Suppose $\eta \le\frac{1}{\bC\kappa^2(f)}$ and let
$\mcX$ be the set constructed by the algorithm for this $\eta$.
Then, for all $x\in \mcG_{\eta}$ either $x\in \mcX$ or
$\|f(x)\| > \delta(f,\eta)$.
\end{lemma}

\proof
Let $x\in\mcG_\eta$. By the definition of $\kappa(f)$
in~\eqref{eq:kappa},
$$
  \frac{1}{\kappa^{2}(f)}\leq 2\max
  \Big\{\mun^{-2}(f,x),\frac{\|f(x)\|^2}{\|f\|^2}\Big\}.
$$
We accordingly divide the proof into two cases.
\medskip

Assume first that
$\max
  \Big\{\mun^{-2}(f,x),\frac{\|f(x)\|^2}{\|f\|^2}\Big\}
=\frac{\|f(x)\|^2}{\|f\|^2}$.

\noindent
In this case
$$
  \eta\le \frac{1}{\bC\kappa^2(f)} \le
  \frac{2\|f(x)\|^2}{\bC\|f\|^2},
$$
which implies
$$
  \|f(x)\|\ge\frac{\sqrt{\eta\,\bC}\,\|f\|}{\sqrt{2}}
  > \frac{\eta\sqrt{\bC}\|f\|}{\sqrt{2}} \geq 1.1\sqrt{(n+1)D}\,\|f\|\,\eta
 = \delta(f,\eta),
$$
the second inequality since $\eta<1$ and the third since
$\bC\ge 12(n+1)D$.
\medskip

Now assume instead that
$\max\Big\{\mun^{-2}(f,x),\frac{\|f(x)\|^2}{\|f\|^2}\Big\}
=\mun^{-2}(f,x)$.

In this case
\begin{equation}\label{eq:mu-eta}
  \eta\le \frac{1}{\bC\kappa^2(f)} \leq
  \frac{2}{\bC\mun^2(f,x)}.
\end{equation}
We will show that the condition
$\frac{1}{531\,\bg(f,x)}\geq\sqrt{\sep(\eta)}$ of the algorithm
holds true,  and that when any of the other two conditions
doesn't hold, then $\|f(x)\|>\delta(f,\eta)$.

Indeed,
$$
  \bg(f,x) = \frac12 D^{3/2}\mun(f,x)
  \underset{\eqref{eq:mu-eta}}{\leq}
  \frac{\sqrt{2}}{2}D^{3/2}\frac{1}{\sqrt{\bC\eta}}\leq
  \frac{1}{531\,\sqrt\eta\,(n+1)^{1/4}}
= \frac{1}{531\,\sqrt{\sep(\eta)}},
$$
the second inequality since
$\sqrt{\bC}\geq \frac{\sqrt{2}}{2}531(n+1)^{1/4}D^{3/2}$.

Assume now that $\ba(f,x)>\alpha_0$. Then
$$
 \a_0 < \frac12 D^{3/2}\mun^2(f,x)\frac{\|f(x)\|}{\|f\|}
$$
which implies
$$
  \|f(x)\| > \|f\| \frac{2\a_0}{D^{3/2}\mun^2(f,x)}
   \underset{\eqref{eq:mu-eta}}{\geq}
   \|f\| \bC\eta\frac{\a_0}{D^{3/2}}
   \geq 1.1\,\sqrt{D(n+1)}\,\|f\| \,\eta=\delta(f,\eta),
$$
the last inequality since
$\bC\geq \frac{531^2}{2} \sqrt{n+1}\,D^{3}
\geq \frac{1.1\sqrt{n+1}D^2}{\a_0}$.

Assume finally that $2.2\,\bb(f,x)\ge\sqrt{\sep(\eta)}$, i.e.
  $$ 2.2 \frac{\|f(x)\|}{\|f\|}\mun(f,x)
      \ge\sqrt\eta(n+1)^{1/4}.
$$
This implies
$$
 \|f(x)\|\ge \|f\|\sqrt \eta\frac{(n+1)^{1/4}}{2.2\,\mun(f,x)}
  \underset{\eqref{eq:mu-eta}}{\geq}
  \|f\|\eta\frac{\sqrt{\bC}(n+1)^{1/4}}{2.2\sqrt{2}}
 \ge 1.1\,\sqrt{(n+1)D}\,\|f\| \,\eta=\delta(f,\eta),
$$
since $\bC\ge  12(n+1)D$.
\eproof

\proofof{Proposition~\ref{prop:algo-covering}}
Lemmas~\ref{lem:density}, \ref{lem:tau-r} and~Remark~\ref{rem:tau} show that
if the algorithm halts, then the current value of
$r$ when halting and that of $\bt:=6r$ satisfy the hypothesis of
Proposition~\ref{prop:SNW}. The fact that $\bt=6r$
 shows that with the choice
$\e:=3.5 r$ the manifold
$\mcM_{\IS}$ is a deformation retract of $U_\e(\mcX)$ and, hence,
the two are homotopically equivalent. Finally,
the fact that $\mcX$ is closed
under the involution $x\mapsto -x$ is straightforward.
This shows correctness.

To evaluate the complexity, note that
Lemma~\ref{lem:L2} shows that the algorithm
halts as soon as
\begin{equation*}\label{eq:eta}
    \eta\leq\eta_0:=\frac{1}{\bC\,\kappa^2(f)}.
\end{equation*}
This gives a bound of $\Oh(\log_2(\bC\kappa(f)))$ for
the number of iterations.

At each such iteration there are at most
$R_{\eta}:=2(n+1)\big(\frac2\eta\big)^n$ points in the grid
$\mcG_{\eta}$. For each such point $x$ we can evaluate
$\mun(f,x)$ and $\|f(x)\|$, both with cost $\Oh(N)$
(cf.~\cite[Prop.~16.45 and Lem.~16.31]{Condition}).
It follows that the cost of each iteration is
$\Oh(R_{\eta}N)$.

Since at these iterations $\eta\geq \eta_0$, we have
$R_{\eta}\leq 2(n+1)\big(2\bC\kappa^2(f)\big)^{n}$. Using this
estimate in the $\Oh(R_{\eta}N)$ cost of each iteration
and multiplying by the bound
$\Oh(\log_2(\bC\kappa(f)))$ for the number of iterations,
we obtain a bound of $N(nD\kappa(f))^{\Oh(n)}$
for the total cost. The claimed bound
follows by noting that $N=(nD)^{\Oh(n)}$.

Finally, the number of points $K$ of the returned $\mcX$ satisfies
\begin{equation}\tag*{\qed}
  K=R_{\eta_0}\le 2(n+1)\big(2 \bC\,\kappa^2(f)\big)^n
 =(nD\kappa(f))^{\Oh(n)}.
\end{equation}

\section{Computing the Betti numbers and torsion coefficients
of spherical and projective algebraic sets}\label{sec:betti}

Let $X$ be a topological space and $\{U_i\}_{i\in I}$ a collection of
open subsets covering $X$.
We recall that the {\em nerve} of this covering
is the abstract simplicial complex $\mcN(U_i)$ defined on
$I$ so that a finite set $J \subset I$ belongs to $\mcN(U_i)$
if and only if  the intersection $\cap_{j\in J} U_j$ is nonempty.
In general the complex does not reflect the topology of $X$, except when intersections are contractible, in which case there is the
Nerve Theorem, that we quote here
from~\cite[Theorem~10.7]{Bjo:95}.

\begin{theorem}\label{thm:nerve}
Let $X$ be a triangulable topological space and $\{U_i\}_{i\in I}$ a
locally finite family of open subsets (or a finite family of closed
subsets) such that $X=\cup_{i\in I}U_i$. If every nonempty finite
intersection $\cap_{j\in J} U_j$ is
contractible, then $X$ and the nerve $\mcN(U_i)$ are
homotopically equivalent.\eproof
\end{theorem}

Here we use the Nerve Theorem in
the case where the sets $U_i$ in the statement of the theorem
are the open balls $B(x_i,\e)$ for $x_i\in \mcX$ where $\{\mcX,\e\}$ is
the output of Algorithm~\ref{alg:covering} and $X$ is their union.
Note that as balls are convex, so is their intersection. Hence,
these intersections, if nonempty, are contractible, and we can apply
the Nerve Theorem.
That is, given $\{\mcX,\e\}$ we want to compute first
 its nerve $\mcN:=\mcN(U_i)$ and then, the Betti numbers and
torsion coefficients of $\mcN$.
Proposition~\ref{prop:algo-covering} and
Theorem~\ref{thm:nerve} ensure that these quantities coincide
for $\mcN$ and $\mcM_{\IS}$.

In what follows, we assume that we have ordered the set $\mcX$
so that $\mcX=\{x_1<x_2<\ldots<x_K\}$ where $K=|\mcX|$ is the
cardinality of $\mcX$. Then,
for $k\geq0$, the abelian group $C_k$ of $k$-chains of
$\mcN$ is free, generated by the set of $k$-faces

\begin{equation}\label{eq:simplices}
   \big\{J\subset \{x_1,\ldots,x_K\}\mid |J|=k
   \mbox{ and } \bigcap_{x_j\in J}B(x_{j},\e)\neq\emptyset\big\}.
\end{equation}

To determine the faces of $C_k$ from $\{\mcX,\e\}$ we need to be able
to decide whether, given a subset $\{x_{i_1},\ldots,x_{i_k}\}$ of
$\mcX$, the intersection of the balls $B(x_{i_j},\e)$,
$j=1,\ldots,k$, is nonempty. This is equivalent to say that the
smallest ball containing all the points $\{x_{i_1},\ldots,x_{i_k}\}$
has radius smaller than $\e$, and we can do so if we have at hand
an algorithm computing this smallest ball.
Since we are looking here for a deterministic algorithm, we do not
apply the efficient but randomized algorithm
of~\cite[pp.~60--61]{EdelHar:10},  whose (expected) cost
is bounded by $\Oh((n+2)!k)$, but  we apply a deterministic
quantifier elimination algorithm to the following problem:
given
$x_{i_1},\ldots,x_{i_k}\in\R^{n+1}$ and $\e>0$, decide whether
$$
   \exists \, z\in\R^{n+1} \mbox{ s.t. }  \|x_{i_j}-z\|<\e
   \mbox{ for } 1\le j\le k.
$$
This can be solved using for instance~\cite{Ren92a}
in  time linear in $k^{\Oh(n)}$.
As there are $\binom{K}{k}\le K^k$
subsets of $k$ elements in $I$, the following result is clear.

\begin{lemma}\label{lem:Ck}
The cost of constructing $C_k$ is bounded by
$K^k \cdot k^{\Oh(n)}$. \eproof
\end{lemma}

For $k\ge 1$ the boundary map $\partial_k:C_k\to C_{k-1}$
is defined, for a
simplex $J\in C_k$, $J=\{x_{i_1},\ldots,x_{i_k}\}$, with
$i_1<i_2<\ldots<i_k$, by
$$
       \partial_k(J)=\sum_{j=1}^k (-1)^j
        \{x_{i_1},\ldots,\widehat{x_{i_j}},\ldots,x_{i_k}\}
$$
where the $(k-1)$-face $ \{x_{i_1},\ldots,\widehat{x_{i_j}},\ldots,x_{i_k}\}$
is obtained by deleting the $j$th element in $J$. This map is
therefore represented
by a matrix $M_k$ with $O_{k-1}$ rows and $O_k$ columns with
entries in $\{-1,0,1\}$,
where $O_k$ denotes the number of faces in~\eqref{eq:simplices}.

\begin{proposition}\label{prop:betti-S}
We can compute the Betti numbers
$b_0(\mcM_\IS),\ldots,b_{n-m}(\mcM_\IS)$
as well as the torsion coefficients of $\mcM_\IS$ with cost
$$
 (n D \kappa(f))^{\Oh(n^2)}.
$$
\end{proposition}

\proof
Algorithm {\sf Covering} produces, as shown in
Proposition~\ref{prop:algo-covering}, a pair $\{\mcX,\e\}$ such that
the union $U_\e(\mcX)$ of the balls $B(x,\e)$, for $x\in\mcX$, covers
$\mcM_\IS$ and is homotopically equivalent to it.
Theorem~\ref{thm:nerve} then ensures that the nerve $\mcN$ of this
covering is homotopically equivalent to $U_\e(\mcX)$ (and hence to
$\mcM_\IS$). It is therefore enough to compute the Betti numbers and
torsion coefficients of $\mcN$.  To do so, we construct, for
$k=0,\ldots, n-m+1$, the group $C_k$ (i.e., we determine its
faces). This has cost
$$
   \sum_{k=0}^{n-m+1} K^k \cdot k^{\Oh(n)}=
   \sum_{k=0}^{n-m+1} (nD\kappa(f))^{\Oh(nk)}
   k^{\Oh(n)}= (n D \kappa(f))^{\Oh(n^2)}
$$
by Lemma~\ref{lem:Ck} and the bound for $K$ in
Proposition~\ref{prop:algo-covering}.

With the groups $C_k$ at hand we write down the matrices $M_k$
corresponding to the boundary maps $\partial_k$, for
$k=1,\ldots,n-m+1$. Next we compute their Smith normal forms
$D_k$,
$$
   D_k=\left[ \begin{array}{cccccc} b_{k,1} & &&&& \\
  &\ddots&&&& \\ && b_{k,t_k} &&&\\ &&& 0 &&\\
  &&&&\ddots&\\ &&&&& 0
   \end{array}\right].
$$
Then, $\dim \Im \partial_k=\Rk(D_k) =t_k$, and consequently $\dim \ker
\partial_k= O_k-\Rk(D_k) =O_k-t_k$.  For $k=1,\ldots,n-m$ we
thus obtain the Betti numbers
$$
  b_k(\mcM_\IS)= \dim \big(\ker
   \partial_k/\Im \partial_{k+1}\big) = O_k-t_k -t_{k+1}
$$
and the same formula yields $b_0(\mcM_\IS)$ and
$b_{n-m}(\mcM_\IS)$ by
taking $t_0=0$.
Furthermore, it is well-known that the $k$th homology group of $\mcN$
(and hence that of $\mcM_\IS$ as well) has the structure
$$
H_k(\mcM_\IS)\simeq\Z^{b_k(\mcM_\IS)} \oplus
\Z_{b_{k+1,1}}\oplus\Z_{b_{k+1,2}}\oplus\ldots
\oplus\Z_{b_{k+1,t_{k+1}}},
$$
that is, its torsion coefficients are
$b_{k+1,1},b_{k+1,2},\ldots,b_{k+1,t_{k+1}}$.

The cost of this last computations
is that of computing the Smith normal forms $D_1,\ldots,D_{n-m}$. The
one for $D_k$ can be done (see~\cite{stor:96}) with cost $$
\Oh{\,\tilde{\ }}\big((\min\{O_k,O_{k-1}\})^5\max\{O_k,O_{k-1}\}
\big)= \Oh{\,\tilde{\ }}\big(K^{6n}\big)= (n D \kappa(f))^{\Oh(n^2)}
$$ (here $\Oh{\,\tilde{\ }}(g)$ denotes $\Oh(g\log^c g)$ for some
constant $c$) and hence the same bound holds for the cost of computing
all of them.
\eproof

The reasoning above extends in a simple manner to compute the
homology
of $\mcM_\P$. Indeed, projective space $\P^n$ is homeomorphic to the
quotient $\IS^n/\sim$ where $\sim$ is the equivalence relation that
identifies antipodal points. Now consider the map
$$
      \IS^n \overset{[~]}{\longrightarrow} \P^n
$$
associating to $x$ its class $[x]=\{x,-x\}$. Because the set $\mcX$
is closed by taking antipodal points, its image $\bar{\mcX}$
under $[~]$ is well-defined and so is the ball in projective space   $B_\P([x],\e):=\{B(x,\e), B(-x,\e)\}$.
Then, the retraction from the union of the balls $B(x,\e)$
onto $\mcM_\IS$ induces a retraction in projective space from the union
of the balls $B_\P([x],\e)$ onto $\mcM_\P$.

Also, given $x_{i_1},\ldots,x_{i_k}$ in $\mcX$, the intersection of
$B([x_{i_j}],\e)$ is nonempty if and only if there exist representatives of
$[x_{i_1}],\ldots,[x_{i_k}]$ such that the Euclidean balls
centered at these representatives
have nonempty intersection. That is, if and only if there exist
$e_1,\ldots,e_k\in\{-1,1\}$ such that the balls
$B(e_1x_{i_1},\e),B(e_2x_{i_2},\e),\ldots,B(e_kx_{i_k},\e)$ have
nonempty intersection. This can be checked
by brute force,  by  checking each of the $2^k$ possibilities.
Furthermore, if this is the case we get, since $\e<1$,
\begin{align*}
   \bigcap_{1\le j\le k}B_\P([x_{i_j}],\e)&= \left[ B(e_1x_{i_1},\e)\cap\ldots \cap B(e_kx_{i_k},\e)\right]\\ &=
   \big\{B(e_1x_{i_1},\e)\cap\ldots \cap B(e_kx_{i_k},\e),
   B(-e_1x_{i_1},\e) \cap \ldots \cap B(-e_kx_{i_k},\e) \big\}.
\end{align*}
Since if $B(e_1x_{i_1},\e)\cap\ldots \cap B(e_kx_{i_k},\e)$ contracts
to $y\in \R^{n+1}$ then
$B(-e_1x_{i_1},\e) \cap \ldots \cap B(-e_kx_{i_k},\e)$ contracts to
$-y$, then  the intersection of
$B([x],\e)$ contracts to $\{y,-y\}=[y]\in \P^n$ and the Nerve Theorem
applies: it implies that the nerve
$\bar{\mcN}$ of the family $\{B([x],\e)\mid [x]\in\bar{\mcX}\}$ is
homotopically equivalent to the union of this family. The reasoning of
Proposition~\ref{prop:betti-S} straightforwardly applies to prove
the following result.

\begin{proposition}\label{prop:betti-P}
We can compute the Betti numbers
$b_0(\mcM_\P),\ldots,b_{n-m}(\mcM_\P)$
as well as the torsion coefficients of $\mcM_\P$ with cost
\begin{equation}\tag*{\qed}
 (n D \kappa(f))^{\Oh(n^2)}.
\end{equation}
\end{proposition}

\section{On the cost of computing coverings for random systems}
\label{sec:random}

The following result is a part of Theorem~21.1
in~\cite{Condition}.

\begin{theorem}\label{th:BCL}
Let~$\Sigma\subset \R^{p+1}$ be contained in a
real algebraic hypersurface, given as
the zero set of a homogeneous polynomial of degree~$d$
and, for $a\in\R^{p+1}$, $a\neq 0$,
$$
    \scC(a):=\frac{\|a\|}{\dist(a,\Sigma)}.
$$
Then, for all $t\ge (2d+1)p$,
\[
   \Prob_{a\in \IS^p}\{\scC(a)\geq t\}\ \leq\
   4e\, dp\, \frac{1}{t}
\]
and
\begin{equation}\tag*{\qed}
  \Exp_{a\in\IS^p} \big(\log_2 \scC(a)\big)\leq
  \log_2 p +\log_2 d + \log_2(4e^2).
\end{equation}
\end{theorem}

\begin{remark}
For condition numbers over the complex numbers, one can
improve the tail estimate in Theorem~\ref{th:BCL} to show a
rate of decay of the order of $t^{-2(p+1-\ell)}$ where $\ell$ is
the (complex) dimension of $\Sigma\subset\C^{p+1}$
(see~\cite[Theorem~4.1]{Demmel88}). Over the reals, such an
estimate (with the 2 in the exponent removed) has only been
proved in the case where $\Sigma$ is
complete intersection~\cite{Lotz15}. We suspect that a similar
estimate holds for $\kappa(f)$.
\end{remark}

We define
$$
   \Sigma_\C:=\Big\{ f\in \Hdm  \,|\, \exists \,x\in  \C^{n+1}
   \mbox{ such that }\sum_{0\le j\le n}x_j^2=1, \ f(x)=0
   \mbox{ and } \Rk(Df(x))<m\Big\}.
$$
The discriminant variety $\Sigma_\R$ defined in~\eqref{def:sigma}
is contained in $ \Sigma_\C$.

\begin{proposition} \label{prop:discrim}
Let $U$ be a set of $N=\dim_\R \Hdm $ variables. Then there exists a
polynomial $G\in \Q[U]\setminus \{0\}$ such that $G|_{\Sigma_\C}=0$
and $\deg(G)\le m^{n+2}(n+1)D^{n+1} $. (Here $G(f)$ for $f\in
\Sigma_\C$ means specializing $G$ at the coefficients of the
polynomials in $f$.)
\end{proposition}

\proof
Observe that for generic $f=(f_1,\dots, f_m)\in \Hdm$ the map
$x\mapsto Df(x)$, $x\in \C^{n+1}$, is surjective, that is
$\Rk(Df(x))=m$, and that the condition $\Rk(Df(x))<m $ is equivalent
to the vanishing of all maximal minors of the matrix $Df(x)\in
\C^{m\times(n+1)}$.

For convenience, we write
$U=\{u_{i,\alpha}\mid i=1,\ldots,m, |\alpha|=d_i\}$. We
consider the general $(n+1)$-variate polynomials of
degree $d_i$,
$$
  F_i=\sum_{|\alpha|=d_i} u_{i,\alpha} X^\alpha \ \in
  \Q[U][X], \quad 1\le i\le m.
$$
Let $DF(U,X)  \in \Q[U][X]^{m\times (n+1)}$ be the Jacobian matrix of
$F=(F_1,\dots,F_m)$ w.r.t.~$X$, and  denote by
$M_k(U,X),\, 1\le k\le t$,  all
its maximal minors.  We consider the polynomials
\begin{equation}\label{pols}
   \sum_{0\le j\le m} X_j^2-1,\, F_i(u_i,X), \,  M_k (U,X), \quad
  1\le i\le m,\,  1\le k\le t.
\end{equation}
These polynomials have no common zeros in
$\overline{\Q(U)}^{n+1}$ because they have no common
zeros for a generic specialization of $U$ as mentioned at the
beginning of the proof, and
we can apply \cite[Cor.4.20]{DKS:13}.  We have
\begin{align*}
   &\deg_X(F_i)=d_i\le D, \, \deg_X\Big(\sum X_j^2-1\Big)=2,\,
   \deg_X(M_k)\le m(D-1),\\
   & \deg_{U}(F_i)=1,\, \deg_{U}\Big(\sum X_j^2-1\Big)=0, \,
   \deg_{U}(M_k)\le m,
\end{align*}
and therefore there exists $G\in \Q[U]\setminus\{0\}$ such that $G$
belongs to the ideal in $\Q[U,X]$ generated by the polynomials
in~\eqref{pols} with
$$
   \deg_{U}(G)\le (mD)^{n+1}\sum_{0\le \ell\le n}m
   \le m^{n+2}(n+1)D^{n+1}.
$$
Clearly this polynomial $G$ vanishes on
all $f\in \Sigma_\C$.
\eproof

\begin{corollary}\label{corol:probs}
Let $\cost_\IS(f)$ and $\cost_\P(f)$
denote the costs of computing the Betti numbers and torsion coefficients
of $\mcM_\IS$ and $\mcM_\P$, respectively. For $f$ drawn from the uniform
distribution on $\IS(\Hd[m])=\IS^{N-1}$ we have the following:
\begin{description}
 \item[(i)]
 With probability at least $1-(nD)^{-n}$ we have
 $\cost_\IS(f)\leq (nD)^{\Oh(n^3)}$.
 Similarly for $\cost_\P(f)$.
 \item[(ii)]
 With probability at least $1-2^{-N}$ we have
 $\cost_\IS(f)\leq 2^{\Oh(N^2)}$.
 Similarly for $\cost_\P(f)$.
\end{description}
\end{corollary}

\proof
For all $t\ge\big(2(n+1)m^{n+2}D^{n+1}+1\big)N$,
it follows from Theorem~\ref{th:BCL} and
Propositions~\ref{prop:kappa-dist}
and~\ref{prop:discrim}, that we have
$$
   \Prob_{f\in \IS^{N-1}}\{\kappa(f)\geq t\}\ \leq\
   4e\, m^{n+2}(n+1)D^{n+1}N\, \frac{1}{t}.
$$
By taking $t=(nD)^{cn}$ for a constant $c$ large enough, we have
$$
   \Prob_{f\in \IS^{N-1}}\{\kappa(f)\geq (nD)^{cn}\}\ \leq\
   4e\, m^{n+2}(n+1)D^{n+1}N\,(nD)^{-cn}\leq (nD)^{-n}.
$$
By Propositions~\ref{prop:betti-S}
and~\ref{prop:betti-P},
for $f$ with $\kappa(f)\leq (nD)^{cn}$ we have
$\cost_\IS(f),\cost_\P(f)\leq (nD)^{\Oh(n^3)}$. This proves~(i).

To prove part~(ii) we take $t=2^{cN}$ for $c$ large enough. Then,
\begin{equation*}
   \Prob_{f\in \IS^{N-1}}\{\kappa(f)\geq 2^{cN}\}\ \leq\
   4e\, (n+1)m^{n+2}D^{n+1}N2^{-cN}\leq 2^{-N}.
\end{equation*}
Using Propositions~\ref{prop:betti-S} and~\ref{prop:betti-P}
again, we have that for $f$ such that $\kappa(f)\leq 2^{cN}$,
$\cost_{\IS}(f),\cost_\P(f)\leq (nD)^{\Oh(n^2)} 2^{\Oh(n^2N)}
\leq 2^{\Oh(N^2)}$, the last as $N\geq \frac{n^2}{2}$.
\eproof

\section{Remaining proofs}\label{sec:proofs}

\subsection{Proof of Proposition~\ref{prop:kappa-dist}}
\label{sec:kappa-dist}

We start by defining a fiber version of $\Sigma_{\R}$.
For $x\in\IS^n$ we let
$$
  \Sigma_\R(x):=\big\{g\in \Hdm: g(x)=0 \mbox{ and }
  \rank(Dg(x))<m\big\}.
$$
Note that, for all $x\in\IS^n$, $\Sigma_\R(x)$ is a
cone in $\R^N$. In particular, $0\in\Sigma_\R(x)$.
The following result is the heart of our proof.

\begin{proposition}\label{prop:kappa-dist-fiber}
For all $f\in\Hdm$ and $x\in \IS^n$,
$$
 \frac{\|f\|}{{\sqrt{2}}\,\dist(f,\Sigma_\R(x))} \le
 \kappa(f,x)\le \frac{\|f\|}{\dist(f,\Sigma_\R(x))}.
$$
\end{proposition}

\proof
{We only need to prove the statement for $f\notin \Sigma_\R(x)$.}
As we saw in Remark~\ref{rem:inv},  $\kappa(\lambda f,x)=\kappa(f,x)$ for all $\lambda\neq 0$, and
also $\dist(\lambda f,\Sigma_\R(x))=|\lambda|\dist(f,\Sigma_\R(x))$.
We can therefore assume, without loss of
generality, that $\|f\|=1$.

Because the orthogonal  group $\scO(n+1)$ in $n+1$ variables
acts on $\Hdm\times\IS^n$
and leaves $\mun,\kappa$ and the distance to $\Sigma_\R$
invariant,
we may assume without loss of generality that $x=e_0:=(1,0,\dots,0)$.

For $1\le i\leq m$ write
\begin{align}\label{fi}
   f_i(X)&=\sum_{q=0}^{d_i} X_0^{d_i-q}f_{i,q}(X_1,\ldots,X_n)
   = X_0^{d_i}f_{i,0} + \sum_{q=1}^{d_i} X_0^{d_i-q}f_{i,q}(X_1,\ldots,X_n) \nonumber \\
   & = X_0^{d_i}f_i(e_0)+ X_0^{d_i-1}\sum_{1\le j\le n}\frac{\partial f_i}{\partial X_j}(e_0)X_j + Q_i(X)
\end{align}
where in the first line $f_{i,q}$ is a homogeneous polynomial of degree $q$, and in the second,
$\deg_{X_0}(Q_i)\le d_i-2$. In particular $f_{i,1}= \sum_{1\le j\le n}\frac{\partial f_i}{\partial X_j}(e_0)X_j$.

\medskip
We first prove that $
   \kappa(f,e_0)\le  {1}/{\dist(f,\Sigma_\R(e_0))}
$, or equivalently,  $$\dist(f,\Sigma_\R(e_0))^2\leq \kappa(f,e_0)^{-2}=\mun^{-2}(f,e_0)+\|f(e_0)\|^2.$$
Write
$f_{i,1}(X_1,\dots,X_n)=\sqrt{d_i}a_{i1}X_1+\cdots+\sqrt{d_i}a_{in}X_n$
for suitable $a_{ij}$. Therefore
$$
  \frac{\partial f_i}{\partial X_j}(e_0) =
  \left\{  \begin{array}{ll}
     d_if_i(e_0)
    & \mbox{if $j=0$}\\[3pt]
 \sqrt{d_i}a_{ij}& \mbox{if $j\geq 1$.}
 \end{array}\right.
$$
Define
$$
  h_i:=f_i-X_0^{d_i}f_{i,0}= \sum_{q=1}^{d_i} X_0^{d_i-q}f_{i,q}(X_1,\ldots,X_n)
$$
for $1\le i\le m$. Then
\begin{equation}\label{h}
 \|f-h\|^2= \sum_{i\leq m} f_{i,0}^2
   = \sum_{i\leq m} f_i(e_0)^2 =\|f(e_0)\|^2.
\end{equation}
In addition $h_i(e_0)=0$ and
 for $0\le j\le n$,
$$
  \frac{\partial h_i}{\partial X_j}(e_0) =
  \left\{  \begin{array}{ll}
    \frac{\partial f_i}{\partial X_j}(e_0) - d_if_i(e_0)\,=\,0
    & \mbox{if $j=0$}\\[3pt]
    \frac{\partial f_i}{\partial X_j}(e_0)=\sqrt{d_i}a_{ij}& \mbox{if $j\geq 1$.}
 \end{array}\right.
$$
Therefore, we have (recall the definition of $\Delta$
from \S\ref{sec:cond-numb})
$$
  \Delta^{-1}Df(e_0) = \left[\begin{array}{cccc}
  \sqrt{d_1}\,f_1(e_0)&a_{11}&\ldots&a_{1n}\\
  \sqrt{d_2}\,f_2(e_0)&a_{21}&\ldots&a_{2n}\\
  \vdots &&& \vdots\\
  \sqrt{d_m}f_m(e_0)&a_{m1}&\ldots&a_{mn}
\end{array}\right]
$$
and
$$
  \Delta^{-1}Dh(e_0) = \left[\begin{array}{cccc}
  0&a_{11}&\ldots&a_{1n}\\
  0&a_{21}&\ldots&a_{2n}\\
  \vdots &&& \vdots\\
  0&a_{m1}&\ldots&a_{mn}
\end{array}\right].
$$
Let $A=(a_{ij})\in\R^{m\times n}$ so that
$\Delta^{-1}Dh(e_0) = [0\, A]$. We know that $\Rk(A)\le m$.
\\
If $\Rk(A)\leq m-1$,  then $h\in \Sigma_\R(e_0)$
and hence, by~\eqref{h}
$$
  \dist(f,\Sigma_\R(e_0))^2\leq \|f-h\|^2 = \|f(e_0)\|^2 \le
\mun^{-2}(f,e_0)+ \|f(e_0)\|^2.
$$
If $\Rk(A)=m$, then (the inequality by~\cite[Lemma~3]{CDW:05}),
\begin{equation}\label{eq:fh}
 \mun(f,e_0)=\|(\Delta^{-1} Df(e_0))^\dagger\|\leq
\|(\Delta^{-1} Dh(e_0))^\dagger\|=\mun(h,e_0).
\end{equation}
Because of the
Condition Number
Theorem~\cite[Corollaries~1.19 and~1.25]{Condition}
there exists a matrix $P\in\R^{m\times n}$ such that
$A+P$ is a non-zero matrix of 
rank less than $m$ and
$$
  \|P\|_F=\|A^\dagger\|^{-1}=\|[0\;A]^\dagger\|^{-1}=
  \|(\Delta^{-1}Dh(e_0))^\dagger\|^{-1}=\mun^{-1}(h,e_0).
$$
Let $E=(e_{ij})=\Delta P\in \R^{m\times n}$ and consider the polynomials
$$g_i(X):=h_i(X) + X_0^{d_i-1}\sum_{j=1}^n e_{ij} X_j, \quad 1\le i\le m.$$
Then $g_i$ are not all zero, $g_i(e_0)=h_i(e_0)=0$,
$\frac{\partial g_i}{\partial X_0}(e_0)
=\frac{\partial h_i}{\partial X_0}(e_0)=0$,
and
$\frac{\partial g_i}{\partial X_j}(e_0)
=\frac{\partial h_i}{\partial X_j}(e_0)+e_{ij}=c_{ij}+e_{ij}=\sqrt{d_i}a_{ij}+e_{ij}$
for $1\le j\le n$.
It follows that
$$
   Dg(e_0)=[0\; \Delta A+E]=[0\;\Delta (A+ P)]
$$
and therefore  $\Rk (Dg(e_0))
<m$. Hence,
$g\in\Sigma_{\R}(e_0)$. In addition,
$$
 \|g-h\|^2 = \sum_{\substack{1\leq i\leq m\\1\leq j\leq n}}
  {d_i\choose {d_i-1,1}}^{-1} e_{ij}^2
  = \sum_{\substack{1\leq i\leq m\\1\leq j\leq n}} d_i^{-1} e_{ij}^2
  = \sum_{\substack{1\leq i\leq m\\1\leq j\leq n}} p_{ij}^2
  = \|P\|_F^2 =\mun^2(h,e_0).
$$
We conclude as
$$
 \|g-f\|^2 = \|g-h\|^2+\|h-f\|^2
 \underset{\eqref{h}}{\,=\,}  \mun^{-2}(h,e_0)+ \|f(e_0)\|^2
 \underset{\eqref{eq:fh}}{\,\leq\,} \mun^{-2}(f,e_0)+ \|f(e_0)\|^2,
$$
and hence, $\dist(f,\Sigma_\R(e_0))^2\leq \|f-g\| ^2
\leq\mun^{-2}(f,e_0)+\|f(e_0)\|^2$.

\medskip

We now prove that
$\kappa(f,e_0)\ge \frac{1}{{\sqrt{2}}\,\dist(f,\Sigma_\R(e_0))}$,
or equivalently, that
$$
   {2}\,\dist(f,\Sigma_\R(e_0))^2\geq \mun^{-2}(f,e_0)+\|f(e_0)\|^2.
$$
Let $g\in \Sigma_\R(e_0)$ be such that
$\dist(f,\Sigma_\R(e_0))^2=\|f-g\|^2$.
As in Identity~\eqref{fi}, write
\begin{align*}
   g_i(X)
   & =  X_0^{d_i-1}\sum_{1\le j\le n}\frac{\partial g_i}{\partial X_j}
   (e_0)X_j + \tilde Q_i(X),
\end{align*}
where we used that $g(e_0)=0$. From this equality
and~\eqref{fi} it follows that
$$
  f_i-g_i = X_0^{d_i} f_i(e_0)
+  \left[X_0^{d_i-1}
  \sum_{1\le j\le n}\left(\frac{\partial f_i}{\partial X_j}(e_0)X_j
- \frac{\partial g_i}{\partial X_j}(e_0)X_j \right)\right]
+ \left[Q_i(X) -\tilde Q_i(X)\right].
$$
As the three terms in this sum do not share
monomials,
\begin{align*}\|f_i-g_i\|^2 &\geq f_i(e_0)^2
  + \sum_{1\le j\le n} \left(\frac{\partial f_i}{\partial X_j}(e_0)
  - \frac{\partial g_i}{\partial X_j}(e_0) \right)^2 \\
 &\geq  \frac{1}{2}f_i(e_0)^2 + \frac{1}{2}\left[\frac{1}{d_i} f_i(e_0)^2 +
   \frac{1}{d_i}\sum_{1\le j\le n} \left(\frac{\partial f_i}{\partial X_j}(e_0)
  - \frac{\partial g_i}{\partial X_j}(e_0) \right)^2\right]
\end{align*}
and hence,
$$
 \|f-g\|^2\geq  \frac{1}{2}\left(\|f(e_0)\|^2 + \left\|
   \mbox{diag}\Big(\frac{1}{\sqrt{d_i}}\Big)Df(e_0)
   - \mbox{diag}\Big(\frac{1}{\sqrt{d_i}}\Big)Dg(e_0)\right\|^2_F\right).
$$

But $\rank\big(\mbox{diag}\big(\frac{1}{\sqrt{d_i}}\big)Dg(e_0)\big)<m$,
and therefore, by the Eckart-Young theorem,
$$
   \left\| \mbox{diag}\Big(\frac{1}{\sqrt{d_i}}\Big)Df(e_0)
  - \mbox{diag}\Big(\frac{1}{\sqrt{d_i}}\Big)Dg(e_0)\right\|_F
  \ge \sigma_m,
$$
the smallest singular value of
$\mbox{diag}\big(\frac{1}{\sqrt{d_i}}\big)Df(e_0)$.
On the other hand,
$$
  \mun(f,e_0)^{-2} =
  \big\|Df(e_0)^{\dagger}\mbox{diag}\big(\sqrt{d_i}\big)\big\|^{-2}
 =\Big\|\Big(\mbox{diag}\Big(\frac{1}{\sqrt{d_i}}\Big)Df(e_0)
   \Big)^\dagger\Big\|^{-2}
 =\Big(\frac{1}{\sigma_m}\Big)^{-2}=\sigma_m^2.
$$
This concludes the proof since
$$
  \|f-g\|^2\ge\frac{1}{2}\big( \|f(e_0)\|^2 +\sigma_m^2\big)
  = \frac{1}{2}\big(\|f(e_0)\|^2 +\mun(f,e_0)^{-2}\big)
$$
as desired.
\eproof

\proofof{Proposition~\ref{prop:kappa-dist}}
We can assume again $\|f\|=1$. We note that
$$
   \dist(f,\Sigma_\R)=\min\{\dist(f,g): g\in \Sigma_\R\}
  =\min\{\dist(f,\Sigma_\R(x)): \, x\in \IS^n\},
$$
since $\Sigma_\R=\bigcup_{x\in \IS^n} \Sigma_\R(x)$.
Then, using Proposition~\ref{prop:kappa-dist-fiber},
\begin{equation*}
  \kappa(f)=\max_{x\in \IS^n} \kappa(f,x)
 \le \max_{x\in \IS^n}\frac{1}{\dist(f,\Sigma_\R(x))}
 =\frac{1}{\displaystyle\min_{x\in \IS^n}\dist(f,\Sigma_\R(x))}
 =\frac{1}{\dist(f,\Sigma_\R)}.
\end{equation*}
Analogously,
\begin{equation}\tag*{\qed}
  \kappa(f)=\max_{x\in \IS^n} \kappa(f,x)
 \ge \max_{x\in \IS^n}\frac{1}{{\sqrt{2}}\,\dist(f,\Sigma_\R(x))}
 =\frac{1}{{\sqrt{2}}\,
   \displaystyle \min_{x\in \IS^n}\dist(f,\Sigma_\R(x))}
 =\frac{1}{{\sqrt{2}}\,\dist(f,\Sigma_\R)}.
\end{equation}

\subsection{Proof of Proposition~\ref{prop:lipschitz2}}
\label{sec:lipschitz2}

The following simple quadratic map, which was introduced by
S.~Smale in  \cite{Smale86},  is useful in several places in
our development,
\begin{equation}\label{eq:def-psi}
  \psi:\ [0,\infty) \,\to \, \R, \quad
         u\,\mapsto \, 1-4u+2u^2.
\end{equation}
It is monotonically decreasing and nonnegative
in $[0,1-\frac{\sqrt 2}{2}]$.

\begin{lemma}\label{lem:lip1}
Let $u:=\|z-y\|\gamma(f,y)$. For all $\e\in(0,1/2]$, if $u\leq \e$ then
$$
   \mun(f,z)\leq \Big(1+\frac52\e\Big)\mun(f,y).
$$
\end{lemma}

\proof
As $Df(y)Df(y)^\dagger=\Id_{\R^m}$ we have
\begin{eqnarray*}
\mun(f,z)&=& \|f\| \|Df(z)^\dagger\Delta\|
\;=\; \|f\| \|Df(z)^\dagger Df(y)Df(y)^\dagger\Delta\| \\
&\le& \|f\| \|Df(z)^\dagger Df(y)\| \|Df(y)^\dagger\Delta\|
\;\le\; \frac{(1-u)^2}{\psi(u)}  \mun(f,y)
\end{eqnarray*}
the last inequality by~\cite[Lemma 4.1(11)]{Bez4}. We now use
that
$$
  \frac{(1-u)^2}{\psi(u)} = 1+u\Big(\frac{2-u}{1-4u+2u^2}\Big)
  \leq 1+\frac52 \e
$$
the last as $u\leq \e\le\frac12$ and the fact that
$\frac{2-u}{1-4u+2u^2}\leq \frac52$ in the interval $[0,\frac12]$.
\eproof

\proofof{Proposition~\ref{prop:lipschitz2}}
Because of~\eqref{eq:HDE} we have
$$
     \|y-z\|\leq \frac{2\e}{D^{3/2}\mun(f,y)} = \frac{\e}{\bg(f,y)}
     \leq \frac{\e}{\gamma(f,y)}.
$$
Hence, we can apply Lemma~\ref{lem:lip1} to
deduce the inequality on the right. For the inequality on the
left, assume it does not hold. That is,
$$
   \mun(f,z)<\frac{1}{1+\frac52\e}\,\mun(f,y) <\mun(f,y).
$$
Then, $\|y-z\|\leq \frac{\e}{\mun(f,y)}\leq \frac{\e}{\mun(f,z)}$
and we can use Lemma~\ref{lem:lip1} with the roles of
$y$ and $z$ exchanged to deduce that
$$
  \mun(f,y)\leq \Big(1+\frac52\e\Big)\mun(f,z)
$$
which contradicts our assumption.
\eproof

\subsection{Proof of Theorem~\ref{prop:tau-gamma}}
\label{sec:tau-gamma}

Recall that $\mcZ$ denotes $f^{-1}(0)\subset\R^{n+1}$ and
$\mcM_\IS=\mcZ\cap\IS^n$. The idea of the proof is to show that if
$p,q\in \mcM_\IS$, $p\neq q$, then there are fixed radius balls
around $p$ and
$q$ such that the normals at $p$ and $q$ to $\mcM_\IS$, i.e., the
normal spaces of their tangent spaces at $\mcM_\IS$, do not intersect
in the intersection of the two balls.  Either the two points are so
far that there will be no intersection between the two balls, or there
are close and in that case, $\mcM_\IS$ around $p$ can be described as
an analytic map by the implicit function theorem. This enables us to
analyze the normals at $p$ and $q$ and their possible intersection.

For the rest of this section we fix an arbitrary  point $p\in\mcM_\IS$,
i.e., such that $f(p)=0$ and $\|p\|=1$,
with a full-rank derivative $Df(p)$ and we set
$\gamma_p:=\max\{\gamma(f,p),1\}$.

For any $\e>0$ and any linear subspace $H\subset\R^{n+1}$
we denote by $B_{\e,H}(0)$ the open $\e$-ball in $H$ centered
at $0$ and by $B_{\e,H}(p):=p+B_{\e,H}(0)$ the same but
centered at $p$. In the special case that $H=\R^{n+1}$ we simply
write $B_\e(0)$ and $B_{\e}(p)$.

We recall that, because of Euler's formula,  $p\in \ker Df(p)$. We define
$$
  T:=\langle p\rangle^\perp, \qquad   H_1:= \ker Df(p)\cap T,
  \qquad  H_2:= \ker Df(p)^\perp \subset T,
  \qquad H_3:=H_2 + \langle p\rangle,
$$
and consider the orthogonal projections $\pi_i:\R^{n+1}\to H_i$
for $i=1,2,3$. Note that  $H_1$, $H_2$, $H_3$  are linear spaces of
dimension $n-m$, $m$, and $m+1$  respectively. In addition,
$T=H_1\perp H_2$ and
$\R^{n+1}=H_1 \perp H_3=\ker Df(p) \perp H_2$,
where the symbol $\perp$ denotes orthogonal
direct sum.

\begin{proposition}\label{prop:M2} Define
$c_1={0.024}$. Then
$\mcZ\cap B_{\frac{c_1}{\gamma_p},T}(p)$ is contained in the
graph of a real analytic
map $\omega:B_{\frac{c_1}{\gamma_p},H_1}(p)\to H_2$
satisfying
$\omega(p)=0$,  $\|D\omega(p+x)\|\leq {2.3}\,\|x\|\gamma_p$
and $\|\omega(p+x)\|\le {1.15}\, \|x\|^2\gamma_p $,
for all $x\in B_{\frac{c_1}{\gamma_p},H_1}(0)$.
\end{proposition}

Figure~\ref{fig:fig6} below attempts to summarize the situation
described in  Proposition~\ref{prop:M2}.
\begin{figure}[H]
\begin{center}
  \begin{tikzpicture}[scale=2]
\def\mypoints{0.6}
\draw[dashed](-0.5,1.2) to [out=70, in=180] (1.5,2);
\draw[dashed](-0.5,1.2) to [out=110, in=0] (-1.5,2);
\draw[dashed](1.5,2) to [out=135, in=45] (-1.5,2);
\draw[very thick](-0.45,2.2) to [out=0, in=270] (0.05,2.4) 
to [out=90, in=0] (-0.15,2.5);
\path [fill=white] (-1.3,2) rectangle(-0.7,2.6);
\path (-1,2.25) node{\small $\mcZ\cap B_{\frac{c_1}{\gamma_p},T}(p)$};
\fill (-0.005,2.27) circle(\mypoints pt) node[right]{\small $\,p$};
\draw (0.5,2.7) -- (-0.14,2.15);
\draw (-1.5,2.8) -- (2.5,2.8);
\draw (-2.5,2.1) -- (1.5,2.1);
\draw (2.5,2.8) -- (1.5,2.1);
\draw (-2.5,2.1) -- (-1.5,2.8);
\draw (0,-0.73) node[below right]{\footnotesize $0$} -- (0,1.6);
\draw (0,2.27) --(0,3.2) node[right]{\small $\langle p\rangle$};
\draw (0.5,-0.3) node[below right]{\footnotesize $H_1$} -- (-0.14,-0.85); 
\draw (0.5, -0.3) -- (0.5, 1) node[right]{\footnotesize $\ker Df(p)$};
\draw (-1.5,-0.2) -- (2.5,-0.2);
\draw (-2.5,-0.9) -- (1.5,-0.9);
\draw (2.5,-0.2) -- (1.5,-0.9);
\draw (-2.5,-0.9) -- (-1.5,-0.2);
\path (1.85,-0.35) node{\small $T:=\langle p\rangle^\perp$};
\path (1.85,2.65) node{\small $T+p$};
\draw (-0.9,-0.57) node[above]{\footnotesize $H_2$} -- (0.2,-0.762); 
\draw (-0.9,2.43) -- (0.2,2.238); 
\end{tikzpicture}
  \caption{}\label{fig:fig6}
\end{center}
\end{figure}

\proof
The general idea is to first apply
(and get explicit bounds for) the Implicit Function Theorem
to get a real analytic
map $\omega_0: B_{\frac{c_1}{\gamma_p},\ker Df(p)}(p)\to H_2$
satisfying that
$\mcZ\cap B_{\frac{c_1}{\gamma_p}}(p)$ is contained in the graph
of $\omega_0$ with $\omega_0(p)=0$,
$\|D\omega_0(p+x)\|\leq 2.3\|x\|\gamma_p$
and $\|\omega_0(p+x)\|\le 1.15\gamma_p\|x\|^2 $ for all
$x\in B_{\frac{c_1}{\gamma_p},\ker Df(p)}(0)$. We then restrict
$B_{\frac{c_1}{\gamma_p}}(p)$ to $B_{\frac{c_1}{\gamma_p},T}(p)$
and $\omega_0$ to $H_1\subset \ker Df(p)$ to obtain $\omega$
satisfying all the stated conditions.

\smallskip
The process is involved and we describe it as a sequence of
claims.

\smallskip
\noindent
{\bf Claim~1.\quad} {\sl For  all $z\in \R^{n+1}$ such that
$u=u(z):=\|z\|\gamma(f,p)<1-\frac{\sqrt 2}{2}$ the derivative
$Df(p+z)|_{H_2}$ of $f$ with respect to $H_2$ at $p+z$ is invertible.}

\smallskip
Indeed,
$$
   \|Df(p)|_{H_2}^{-1}Df(p+z)|_{H_2}-\Id_{H_2} \|\le \|Df(p)^\dagger
   Df(p+z)-\pi_2\|< \frac{1}{(1-u)^2}-1<1
$$
the first inequality by properties of Moore-Penrose inverse
and the second by~\cite[Lem.~4.1(9)]{Bez4}.
Therefore, by~\cite[Lem.~15.7]{Condition},
$Df(p)|_{H_2}^{-1}Df(p+z)|_{H_2}$ is invertible, which implies
$Df(p+z)|_{H_2}$  invertible as desired. This proves Claim~1.

\smallskip

From now on, since  $\R^{n+1}= \ker Df(p)\oplus H_2$, we  write
indistinctly $f(p)$ or $f(p,0)$ as $p\in \ker Df(p)$, and for
$z=(x,y)\in \ker Df(p)\oplus H_2$, $f(p+z)$ or $f(p+x,y)$.

Let
$$
  \Omega:=\Big\{z=(x,y)\in\ker Df(p)\oplus H_2 \mid \|z\|\leq
\Big(1-\frac{\sqrt{2}}{2}\Big)\frac{1}{\gamma(f,p)}\Big\}.
$$
For all $z=(x_0,y_0)\in\Omega$,
Claim~1 ensures that $Df(p+x_0,y_0)|_{H_2}$ is invertible.
If $f(p+z)=0$, the {Analytic}
Implicit Function Theorem ensures the existence
of an open set $U\subset\ker Df(p)$ around
$x_0$, an open set $V\subset H_2$ around $y_0$ and a
real analytic map ${\omega}_z: p+U\to V$ such that
\begin{equation}\label{eq:sigma-z}
  \{(p+x,{\omega}_z(p+x))\mid x\in U\}
  =\{(p+x,y)\in (p+U)\times V\mid f(p+x,y)=0\}.
\end{equation}
\smallskip

Recall the decreasing map  $\psi$ defined in~\eqref{eq:def-psi} and
consider also the function
$$
  \phi(u):=\frac{(1-u)^2}{\psi(u)} \Big(\frac{1}{(1-u)^2}-1\Big)
  = \frac{2u(1-\frac{u}2)}{\psi(u)}.
$$
We observe that $\phi(u)<2.2\, u$ for $u<0.024=:c_1$.
\smallskip

\noindent
{\bf Claim~2.\quad}
{\sl Let $z=(x_0,y_0)\in\Omega$ and $u=u(z)=\|z\|\gamma(f,p)$.
If $f(p+z)=0$ then
$\|D{\omega}_z(p+x_0)\|\leq \phi(u)$.}
\smallskip

This is Lemma 5.1 in~\cite{Bez4}
(with $x,y$ and $\sigma$ there corresponding to
$p,p+z$ and $D{\omega}_z$ in our context,
and in the particular case where $f(p+z)=0$).
\smallskip

Let now $\omega_0$ be $\omega_z$ for
$z=(0,0)$, and denote by $0\in U_0\subset \ker Df(p)$ and
$0\in V_0\subset H_2$ the open sets given by the Implicit
Function Theorem in last paragraph.
We observe that by Claim~2, we have $D\omega_0(p)=0$
since  $\phi(0)=0$.
\smallskip

\noindent
{\bf Claim~3.\quad}
{\sl We have
$$
   \|D^2{\omega}_0(p)\| \leq 2\gamma(f,p).
$$}

First note that by the Implicit Function Theorem,
$D{\omega}_0(p)=- (Df(p,0)|_{H_2})^{-1} \circ Df(p,0)|_{\ker Df(p)}=0$
and $f\circ(\Id,{\omega}_0)=0$ in $(p+U_0)\times V_0$, so
\begin{eqnarray*}
0&=& D^2(f\circ(\Id,{\omega}_0))(p)\\
&=& D^2f((\Id,0),(\Id,0))(p)+(Df(\Id,\omega_0)(0,D^2{\omega}_0))(p)\\
&=& D^2f((\Id,0),(\Id,0))(p)+Df(p)|_{H_2} D^2{\omega}_0(p).
\end{eqnarray*}
Note we have removed the symbol $\circ$ in the compositions from
the second line above. We have done, and keep doing, this to
make the notation lighter.

So, $D^2{\omega}_0(p)=-Df(p)^{\dagger}D^2f((\Id,0),(\Id,0))(p)$
and we obtain the inequality
$$
   \|D^2{\omega}_0(p)\|= \|Df(p)^{\dagger}D^2f((\Id,0),(\Id,0))(p)\|
   \leq 2 \max_{k>1}\left\|Df(p)^\dagger \frac{D^kf(p)}{k!}\right\|^{1/k-1}
  = 2\gamma(f,p)
$$
from the definition of $\gamma(f,p)$. Claim~3 is proved.
\smallskip

\noindent
{\bf Claim~4.\quad} {\sl Recall $c_1=0.024$.
The analytic map ${\omega}_0: p+U_0\to V_0$ can be analytically extended
on the open ball $B_{\frac{c_1}{\gamma_p},\ker Df(p)}(p)$, and
for all  $p+x\in B_{\frac{c_1}{\gamma_p},\ker Df(p)}(p)$, its extension
--also denoted by ${\omega}_0$-- satisfies the following:
\begin{description}
\item[(i)]
$\|D\omega_0(p+x)\| < 2.3\,\|x\|\gamma_p$, and
\item[(ii)]
$\|\omega_0(p+x)\| < 1.15\,\|x\|^2\gamma_p
< \frac{0.0007}{\gamma_p}$.
\end{description}}
\smallskip

Since $\omega_0$ is defined in $p+U_0$, there exists $r$, $0<r\le
\frac{c_1}{\gamma_p}$, such that $\omega_0$ is defined on
$B_{r,\ker Df(p)}(p)$ and satisfies Conditions~(i) and (ii).
To see~(i) we note that the equality $\|D\omega_0(p)\|=0$
along with Claim~3, the Mean Value Theorem and the fact that
${\omega}_0$ is defined and $C^2$ on $p+U_0$ imply that
\begin{equation}\label{eq:A}
     \|D\omega_0(p+x)\|< 2.3 \|x\| \gamma_p
\end{equation}
for $x$ sufficiently close to $0$. For (ii), from~\eqref{eq:A} and the
Fundamental Theorem of Calculus, we have
\begin{eqnarray}\label{eq:norm-sigma*}
\|\omega_0(p+x)\|&=&
  \|\omega_0(p+x)-\omega_0(p)\|
  \,\leq  \,\int_0^{1}\|D\omega_0(p+tx)x\|dt\nonumber\\
  & \leq&   \int_0^1\|D\omega_0(p+tx)\| \|x\| dt
     \;\underset{\mathrm{(i)}}{<}\;
     \int_0^1 2.3\,t\,\|x\|^2\gamma_p\, dt\\
 &=& 1.15\,\|x\|^2\gamma_p
  \; \leq \;1.15\,\frac{c_1^2}{\gamma_p}
    \;< \; \frac{0.0007}{\gamma_p}. \nonumber
\end{eqnarray}

Let us show that the supremum $r_0$ of all $0<r\le
\frac{c_1}{\gamma_p}$ such that $\omega_0(p+x)$ can be
analytically extended to
$B_{r,\ker Df(p)}(p)$ and satisfies Conditions~(i) and (ii) is exactly
$r_0=\frac{c_1}{\gamma_p}$.  We assume the contrary, that
$r_0<\frac{c_1}{\gamma_p}$, and show that in that case
$\omega_0$ can be extended a little further.

Let $x_0$ be any point in $\ker Df(p)$ with
$\|x_0\|=r_0$. We note that the continuous map $\omega_0$
is bounded on the ball $B_{r_0,\ker Df(p)}(p)$ because
Condition~(i) holds there. Thus we can consider the limit
$y_0:=\displaystyle\lim_{t\to 1^-}{\omega}_0(p+tx_0)$.
Then, reasoning as in~\eqref{eq:norm-sigma*}
\begin{equation*}
\|y_0\|\leq \int_0^{1}\|D\omega_0(p+tx_0)x_0\|dt
\,<\, 1.15\,\|x_0\|^2\gamma_p
    \,\leq\, 1.15\,\frac{c_1^2}{\gamma_p}
    \,< \, \frac{0.0007}{\gamma_p}.
\end{equation*}
Using the inequality above  and the
triangle inequality we obtain
\begin{eqnarray}\label{eq:norm-z}
    \|(x_0,y_0)\|\gamma_p
   &{<} &\Big(\|x_0\|+1.15\,\|x_0\|^2\gamma_p\Big)\gamma_p
   = \|x_0\|\gamma_p \Big(1+1.15\,\|x_0\|\gamma_p\Big)\nonumber\\
   &\underset{r_0<\frac{c_1}{\gamma_p}}{\le}&
   c_1 \Big(1+1.15\, c_1\Big)
   < \Big(1-\frac{\sqrt{2}}{2}\Big).
\end{eqnarray}
Hence, $z:=(x_0,y_0)\in\Omega$ and
$f(x_0,y_0)=0$. This implies that there exist an open ball
$U\subset \ker Df(p)$ and an open set $V\subset H_2$ around
$x_0$ and $y_0$  respectively, and a real analytic
${\omega}_z:p+U\to V$ such that~\eqref{eq:sigma-z} holds.

Since $\|y_0\|<1.15\,\|x_0\|^2\gamma_p$, by taking a smaller
ball $U$  we can further ensure
that ${\omega}_z(p+x)\subset B_{1.15\,\|x\|^2\gamma_p,H_2}(0)$
for $x\in U$. So,~(ii) holds for ${\omega}_z$ on $p+U$.
Furthermore, we may
use Claim~2, Inequality~\eqref{eq:norm-z}, and the fact that
$\phi(u)<2.2\,u$ for all $0<u<c_1$ to deduce
$$
  \|D\omega_z(p+x_0)\|
  < 2.2\|x_0\|\gamma_p\Big(1+1.15\,c_1\Big)
  < 2.3\,\|x_0\|\gamma_p,
$$
so that  ${\omega}_z$ also satisfies~(i)
on~$p+U$, possibly taking an even smaller $U$.

Finally, since the analytic maps $\omega_0$, defined on $p+x\in
B_{r_0,\ker Df(p)}(p) $, and $\omega_z$, defined on $p+x_0+x$ for
$x-x_0\in U$, coincide by~\eqref{eq:sigma-z} on
$B_{r_0,\ker Df(p)}(p) \cap U$, which is non-empty and connected,
${\omega}_z$ is an analytic continuation of ${\omega}_0$ on
$p+U$ around $p+x_0$.
Let us denote by $U_x$ this open ball around $x$ for
$x\in\IS_0:=\{x\in\ker Df(p)\mid \|x\|=r_0\}$ and let
$\mcU:=\cup_{x\in\IS_0} U_x$. Consider the function
$\varphi:\IS_0\to\R$ defined by
$$
    \varphi(x)=\sup\{t\in \R\,\mid \,[x,tx))\subset \mcU\}.
$$
Note that by construction $t>1$ for every $x$.
As $\varphi$ is continuous and $\IS$ is compact and connected
the image $\varphi(\IS)$ is a closed interval $[\ell,\ell']$ with
$1\leq \ell\leq \ell'$.
Furthermore, there exists $x_*\in\IS_0$ such that
$\varphi(x_*)=\ell$ and, hence, $\ell\ge \frac{r_0+r_*}{r_0}>1$,
where $r_*$ is the radius of $U_{x_*}$. It follows that
we can extend ${\omega}_0$ to the open ball in $\ker Df(p)$
centered at $p$ with radius $r_0+r_*>r_0$ and both (i) and (ii)
hold in this ball, a contradiction.
This finishes the proof of Claim~4.
\smallskip

Claim~4 shows that for all
$x\in B_{\frac{c_1}{\gamma_p},\ker Df(p)}(0)$
the point $y={\omega}_0(p+x)$ satisfies $(p+x,y)\in\mcZ$ and
$\|y\|\leq\frac{0.0007}{\gamma_p}$. We will next see that
it is the only point in $H_2$ satisfying these two conditions.
To do so, for each $x\in \ker Df(p)$, we define
$g_x:H_2\to\R^m$ as the
restriction of $f$ to $\{p+x\}\times H_2$ so that $g_x(0)=f(p+x)$.
Because of Claim~1, for all $y\in H_2$ such that
$\|(x,y)\|<\frac{1-\frac{\sqrt{2}}{2}}{\gamma_p}$,
$Df(p+x,y)\mid_{H_2}$ is  invertible. In particular,
$Dg_x(0)=Df(p+x)|_{H_2}$ is invertible
for $\|x\|<\frac{1-\frac{\sqrt 2}{2}}{\gamma_p}$.
\smallskip

\noindent
{\bf Claim~5.\quad}
{\sl  For all $x\in \ker Df(p)$ such that
$u=u(x)=\|x\|\gamma(f,p)<1-\frac{\sqrt 2}{2}$, we have
$\alpha(g_x,0)\le \frac{u}{\psi(u)^2}$.}

\smallskip

To show this claim we adapt the proof
of~\cite[Prop.~3, p.~160]{BCSS98}. First we verify that
$\gamma(g_0,0)=\gamma(f|_{\{p\}\times H_2},p)\le \gamma(f,p)$.
To do this we note that
$$
   \gamma(f,p):=\max_{k>1}\left\|Df(p)^{\dagger}
   \frac{D^kf(p)}{k!}\right\|^{\frac{1}{k-1}}
  =\max_{k>1}\max_{w_1,\dots,w_k\in\IS^n}
   \left\|Df(p)^\dagger \frac{D^kf(p)}{k!}(w_1,\dots,w_k)\right\|
$$
and
\begin{align*}
  \gamma(f|_{\{p\}\times H_2},p)
  :&=\max_{k>1}\left\|Df|_{\{p\}\times H_2}(p)^{-1}
  \frac{D^kf|_{p\times H_2}(p)}{k!}\right\|^{\frac{1}{k-1}}\
  = \ \max_{k>1}\left\|Df(p)|_{H_2}^{-1}
  \frac{D^kf(p)|_{H_2}}{k!}\right\|^{\frac{1}{k-1}}\\
 &= \max_{k>1}\max_{v_1,\dots,v_k\in \IS^{m-1}}
      \left\|Df(p)|_{H_2}^{-1} \frac{D^kf(p)|_{H_2}}{k!}
       (v_1,\dots,v_k)\right\|^{\frac{1}{k-1}}.
\end{align*}
Modulo an orthogonal change of basis (that does not modify
norms), we can write $Df(p)^{\dagger}
=\left(\begin{array}{c}0\\ Df(p)|_{H_2}^{-1}\end{array}\right)$. This
proves that $\gamma(g_0,0)\le \gamma(f,p)$. Also,
$$
  \beta(g_x,0)= \|Dg_x(0)^{-1}\,g_x(0)\|\le \|Df(p+x)|_{H_2}^{-1}
  \,Df(p)|_{H_2}\|\|Df(p)|_{H_2}^{-1} \,f(p+x)\|.
$$
By~\cite[Lem.~4.1(10)]{Bez4},
$$
  \|Df(p+x)|_{H_2}^{-1} \,Df(p)|_{H_2}\|\le
  \frac{(1-u)^2}{\psi(u)}
$$
while by the multivariate version
of~\cite[Lem.~4(b), p.~161]{BCSS98},
$$
  \|Df(p)|_{H_2}^{-1} \,f(p+x)\|\le \frac{\|x\|}
   {1-\|x\|\gamma(f|_{\{p\}\times H_2}, p) }\le \frac{\|x\|}{1-u},
$$
since $\beta(f,p)=0$. This implies
$\beta(g_x,0) \le\frac{1-u}{\psi(u)}\|x\|$.

Also, in the same way that we verified that
$\gamma(g_0,0)\le \gamma(f,p)$ we can check that
$\gamma(g_x,0)\le \gamma(f,p+x)$, and therefore, as in the proof
of~\cite[Prop.~3, p.~162]{BCSS98}, one gets
\begin{equation}\label{eq:2gammas}
   \gamma(g_x,0)\le\frac{\gamma(f,p)}{\psi(u)(1-u)}.
\end{equation}

Multiplying $\beta(g_x,0)$ and $\gamma(g_x,0)$ we conclude
that as long as $u=\|x\|\gamma(f,p)<1-\frac{\sqrt 2}{2}$ we have
$$
   \alpha(g_x,0)\le\frac{u}{\psi(u)^2}.
$$
This proves Claim~5.
\smallskip

\noindent
{\bf Claim~6.\quad} {\sl Recall $c_1=0.024$. For all $x\in\ker Df(p)$ with
$\|x\|\le\frac{c_1}{\gamma_p}$,
there is at most one zero of
the map $g_x$ in the ball
$\|y\|<\frac{0.044}{\gamma_p}$.}
\smallskip

For $0\le u\le c_1$ one has $0.905 \le \psi(u)\le 1$ and
$\frac{u}{\psi(u)^2}< 0.03$. Thus, by Claim~5,
$\alpha(g_x,0)<0.03$ for all $x\in \ker Df(p)$ with
$\|x\|\le\frac{c_1}{\gamma_p}$. The second statement
in Theorem~\ref{alpha-th} applied to $g_x$
tells us that $0$ converges to a
zero of $g_x$ and that all points in the ball of radius
$\frac{0.05}{\gamma(g_x,0)}$ converge to the same zero.
This implies that there is at
most one zero of $g_x$ in the ball of radius
$$
   \frac{0.05}{\gamma(g_x,0)}
   \underset{\eqref{eq:2gammas}}{\geq}
   \frac{0.05\,\psi(u)(1-u)}{\gamma(f,p)}\geq
   \frac{0.05\,\psi(0.024)(1-0.024)}{\gamma(f,p)}
   \geq  \frac{0.044}{\gamma(f,p)}
$$
which proves Claim~6.
\smallskip

We can now finish the proof of the proposition.
Since $B_{\frac{c_1}{\gamma_p}}(p)\subset
B_{\frac{c_1}{\gamma_p},\ker Df(p)}(p) \times
B_{\frac{0.044}{\gamma_p},H_2}(0)$,
it follows from Claims~4 and~6 that
$\mcZ\cap B_{\frac{c_1}{\gamma_p}}(p)$ is included
in the graph $\Gr(\omega_0)$ of $\omega_0$.
We finally restrict $\mcZ\cap B_{\frac{c_1}{\gamma_p}}(p)$ to
$\mcZ \cap B_{\frac{c_1}{\gamma_p},T}(p)$, and therefore
$\omega_0$ restricts to
$\omega: B_{\frac{c_1}{\gamma_p},H_1}(p)\to H_2$,
as explained at the beginning of the proof. The bounds
for $\|D\omega (p+x)\|$ and $\|\omega(p+x)\|$ follow from Claim~4.
\eproof

\smallskip
%

\begin{lemma}\label{lem:inject}
Let $\omega$ be the map of Proposition~\ref{prop:M2} and define
the following continuous map
\begin{equation*}\label{Phi}
  \Phi: B_{\frac{c_1}{\gamma_p},H_1}(0)\subset {H_1}
  \longrightarrow{H_1}, \quad
  \Phi(x)= \frac{x}{\|(p+x,\omega(p+x))\|}.
\end{equation*}
Then  $\Phi$
is a bijection onto its image and satisfies
\begin{description}
\item[(i)]
$\|\Phi(x)\|\geq 0.9997\|x\|$,
\item[(ii)]
$\|D\Phi(x)^{-1}\|\leq 1.0013$.
\end{description}
\end{lemma}

\proof
If we define the map
\begin{equation}\label{S}
  S:B_{\frac{c_1}{\gamma_p},H_1}(0)\to\R, \
  S(x)=\|(p+x,\omega(p+x))\|,
\end{equation}
then $\Phi(x)=\frac{x}{S(x)}$, which implies that $\Phi$ maps rays to
themselves.  To see that $\Phi$ is bijective, it is therefore
sufficient to see that it is monotone increasing along rays, so we
study its derivative along rays and show it is positive.

Let $x=tv$ with $v$ a unit vector and differentiate
$\Phi(tv)=\frac{tv}{S(tv)}$ w.r.t. $t$ to obtain
$$
 \frac{d\Phi}{dt}(tv)=
 \frac{(1+\|\omega(p+tv)\|^2
   -t\langle D\omega(p+tv)v,\omega(p+tv)\rangle)}{S(tv)^3}v.
$$
As we have
\begin{eqnarray*}
  t|\langle D\omega(p+tv)v,\omega(p+tv)\rangle | &\leq&
  t\|D\omega(p+tv)\| \|\omega(p+tv)\|
  \underset{\mathrm{Prop.~\ref{prop:M2}}}{\leq}
  2.3\,t^2\gamma_p \|\omega(p+tv)\|\\
&\underset{t\le\frac{c_1}{\gamma_p}}{\leq}&
  2.3\,\frac{c_1^2}{\gamma_p} \|\omega(p+tv)\|
 < 2\|\omega(p+tv)\|
\end{eqnarray*}
since $c_1=0.024$ and $\gamma_p\ge 1$,
it follows that
\begin{eqnarray*}
 1+\|\omega(p+tv)\|^2-t\langle D\omega(p+tv)v,\omega(p+tv)\rangle
 &>&
 1+\|\omega(p+tv)\|^2-2\|\omega(p+tv)\|\\
&=& (1-\|\omega(p+tv)\|)^2\geq 0,
\end{eqnarray*}
since $\|\omega(p+tv)\|\le 1.15 t^2\gamma_p \le 1.15 c_1^2<1$
by Proposition~\ref{prop:M2} for $|t|\le \frac{c_1}{\gamma_p}$.
This shows that $\Phi$ restricted to
$\{tv\}_{|t|\leq \frac{c_1}{\gamma_p}}$
is strictly monotone, as  wanted.
\smallskip

To show the bounds (i--ii) we first note that for any $x$ with
$\|x\|<\frac{c_1}{\gamma_p}$, by Proposition~\ref{prop:M2},
we have
\begin{eqnarray*}\label{eq:S(x)}
 S(x)&=&(1+\|x\|^2+\|\omega(p+x)\|^2)^{\frac12}
 \le \sqrt{1+\frac{c_1^2}{\gamma_p^2}+ {1.15^2}
 \frac{c_1^4}{\gamma_p^2}}
 \le \sqrt{1+c_1^2+ {1.15^2}c_1^4} \le 1.0003,
\end{eqnarray*}
and hence $\|\Phi(x)\|=\frac{\|x\|}{S(x)}\geq 0.9997\|x\|$. This shows (i).
\smallskip

Also, for any $y\in H_1$,
$$
  DS(x)\,y= \frac{\langle x,y\rangle+
   \big\langle \omega(p+x),D\omega(p+x)y\big\rangle}{S(x)}.
$$
As $ \langle x,y\rangle \leq \|x\| \|y\|$ and by Proposition~\ref{prop:M2},
$$
  \big\langle \omega(p+x), D\omega(p+x)y \big\rangle\leq
 \|\omega(p+x)\| \big\|D\omega(p+x)\big\|\|y\| \leq
  {2.65}\|x\|^3\gamma_p^2\|y\|
  \underset{\|x\|< \frac{c_1}{\gamma_p}}{\leq} 2.65\,c_1^2 \|x\| \|y\|,
$$
we deduce that
\begin{equation*}
\big\|DS(x)y\big\|
\leq
\frac{1}{S(x)} (1+{2.65}c_1^2)\|x\|\|y\|\le
 \frac{1.0016\|x\|\|y\|}{S(x)} .
\end{equation*}
So,
\begin{equation}\label{DS}
  \|DS(x)\| \leq \frac{1.0016\|x\|}{S(x)}  \le 1.0016\|x\|
\end{equation}
since $S(x)\ge 1$.

We now use that $\Phi(x)=\frac{x}{S(x)}$ to derive that,
for any $y\in H_1$,
$$
  D\Phi(x)y =
   \frac{ S(x)  -x DS(x) }
   {S(x)^2}y
$$
and, consequently,
\begin{eqnarray*}
 \big\|D\Phi(x)y\big\|
 \geq  \frac{S(x)-1.0016\|x\|^2}{S(x)^2}\|y\|
 \underset{{S(x)}\ge 1}{\geq}
 \frac{1-1.0016\,c_1^2}{S(x)^2}\|y\|
 \geq
 \frac{1-1.0016\,c_1^2}{1.0003^2}\|y\|
\end{eqnarray*}
the last inequality since $S(x)\le  1.0003$. Therefore,
$$
  \big\|D\Phi(x)y\big\|\geq 0.9988\|y\|.
$$

It follows that the smallest singular value $\sigma$ of $D\Phi(x) $
satisfies $\sigma\ge 0.9988$ and therefore
\begin{equation*}
\big\|D\Phi(x)^{-1}\big\|= \frac{1}{\sigma} \leq 0.0013.
\end{equation*}
This shows (ii).
\eproof

\smallskip
In what follows, we  denote by $\phi$ the map
$\phi: \R^{n+1}\setminus \{0\}\to \IS_n$, $\phi(z)=\frac{z}{\|z\|}$,
as we did in  Section~\ref{sec:grids}.

\begin{lemma}\label{lem:M1}
For all $\e\in(0,1)$,
$\phi(\mcZ\cap B_{\e,T}(p))=\mcM_\IS\cap\phi(B_{h(\e)}(p))$,
where $h(\e):=\frac{\e}{\sqrt{1+\e^2}}$.
\end{lemma}

\proof
The worst possible situation corresponds to a point
$z\in\mcZ\cap B_{\e,T}(p)$ with $\|z-p\|=\e$. This
situation is depicted in Figure~\ref{fig:fig5}.
\begin{figure}[H]
\begin{center}
  \begin{tikzpicture}[scale=5]
\def\mypoints{0.3}
\fill (0,0) circle(\mypoints pt) node[left]{\small $0$};
\fill (0.355,0.355) circle(\mypoints pt) node[above right]{\small $p$};
\fill (0.185,0.525) circle(\mypoints pt) node[left]{\small $z$};
\draw (-0.5,0) arc(180:0:0.5);
\draw (0.355,0.355) circle(0.215);
\draw (0,0.71) -- (0.71,0);
\draw (0,0) -- (0.23,0.65);
\draw (0,0) -- (0.355,0.355);
\draw[thick] (0.155,0.435) -- (0.355,0.355);
\draw[thick] (0.185,0.525) -- (0.355,0.355);
\path (0.32,0.45) node{\small $\e$};
\path (0.22,0.35) node{\small $h(\e)$};
\end{tikzpicture}
  \caption{}\label{fig:fig5}
\end{center}
\end{figure}
If $\a$ denotes the angle at the origin in the figure,
then $\e=\tan(\a)$ and $h(\e)=\sin(\a)$ so that
$h(\e)=\sin\arctan(\e)=\frac{\e}{\sqrt{1+\e^2}}$.
\eproof

\begin{proposition}\label{prop:M3}
Define $c_2= 0.023$ and let $\Phi$ be the map defined in
Lemma~\ref{lem:inject}. Then
$\mcM_\IS\cap \phi\Big(B_{\frac{c_2}{\gamma_p}}(p)\Big)$
is contained in  the graph of a real analytic map
$$
   \vartheta\colon p+ \Phi\Big(B_{\frac{c_1}{\gamma_p},H_1}(0)\Big)
   \subset p+H_1 \to H_3
$$
satisfying $\vartheta(p)=0$,
$\|D\vartheta(p+x)\|\leq {3.4}\|x\|\gamma_p$ and
$\|\vartheta(p+x)\|\le  {1.7}\|x\|^2\gamma_p$, for all
$x\in \Phi\Big(B_{\frac{c_1}{\gamma_p},H_1}(0)\Big)$.
Moreover,
$B_{\frac{c_2}{\gamma_p},H_1}(0) \subset
\Phi\Big(B_{\frac{c_1}{\gamma_p},H_1}(0)\Big)$.
\end{proposition}

\proof
Write $\mcB= \Phi\Big(B_{\frac{c_1}{\gamma_p},H_1}(0)\Big)$.
By Lemma~\ref{lem:inject}, $\Phi$ is a bijection onto $\mcB$.
Let $\omega$ be the map defined in Proposition~\ref{prop:M2}.
We recall $H_3=H_2+\langle p\rangle$ and $\pi_3$ is the
projection onto $H_3$ and define
\begin{eqnarray*}
  \vartheta:p+\mcB&\to& H_3\\
  p+x&\mapsto& \big(\pi_3 \phi(\Id,\omega)\big)(p+\Phi^{-1}(x))-p,
\end{eqnarray*}
where $\phi(p+x',y):=\frac{(p+x',y)}{\|(p+x',y)\|}$
for $(x',y)\in H_1\times H_2$.\\
Note that $\vartheta(p)=0$.\\
For  $x':=\Phi^{-1}(x)$ we have
$$
  \vartheta(p+x)=\pi_3\phi(p+x',\omega(p+x'))-p.
$$
Also note that
$x=\Phi(x')=\frac{x'}{\|(p+x',\omega(p+x')\|}
=\pi_1\phi(p+x',\omega(p+x'))$ implies
that, for each $x\in\mcB$ (or, equivalently, for
each $x'\in \Phi^{-1}(\mcB)=B_{\frac{c_1}{\gamma_p},H_1}(0)$) ,
\begin{eqnarray}\label{eq:graf}
         \big(p+x,\vartheta(p+x)\big)&=&
         \big(p+x,\pi_3\phi(p+x',\omega(p+x'))-p\big)=
  \big(x,\pi_3\phi(p+x',\omega(p+x'))\big)\nonumber\\
&=&\big(\pi_1,\pi_3)\phi(p+x',\omega(p+x')\big)
                  = \phi\big(p+x',\omega(p+x')\big)
\end{eqnarray}
modulo the identification $H_1\times H_3=H_1\oplus H_3=\R^{n+1}$.
Identity~\eqref{eq:graf} shows that
$\Gr(\vartheta)=\phi(\Id,\omega)(B_{\frac{c_1}{\gamma_p},H_1}(p))$.

Now, from Proposition~\ref{prop:M2} we know that
$$
  \mcZ\cap B_{\frac{c_1}{\gamma_p},T}(p)\subseteq
  \Gr(\omega)=(\Id,\omega) (B_{\frac{c_1}{\gamma_p},H_1}(p))
$$
and therefore, by Lemma \ref{lem:M1},
$$
  \mcM_\IS\cap \phi\Big(B_{\frac{c_1}{\sqrt{\gamma_p^2+c_1^2}}}(p)\Big)
  = \phi(\mcZ\cap B_{\frac{c_1}{\gamma_p},T}(p))\subseteq
  \phi(\Id,\omega)(B_{\frac{c_1}{\gamma_p},H_1}(p))=\Gr(\vartheta).
$$
As $\gamma_p\ge 1$ we have
$\gamma_p^2+c_1^2\le {1.0006} \gamma_p^2$ and
therefore
$ \frac{c_1}{\sqrt{\gamma_p^2+c_1^2}}\ge \frac{c_2}{\gamma_p}$
for $c_2:={0.023}$.
This shows that $\mcM_\IS\cap
\phi\Big(B_{\frac{c_2}{\gamma_p}}(p)\Big) \subset \Gr(\vartheta)$.
\smallskip

We now show the bounds. By definition, for all $x\in \mcB$ one
has $\vartheta(p+x)=(\psi_3\circ \Phi^{-1}) (x)$, where
$\psi_3:B_{\frac{c_1}{\gamma_p},H_1}(0)\to H_3$ is defined as
$$
  \psi_3(x'):=\pi_3\phi(\Id,\omega)(p+x')-p=
   \frac{\omega(p+x')}{S(x')}-\Big(1-\frac{1}{S(x')}\Big)p,
$$
where $S(x')$ is defined in \eqref{S}.
Hence, for $x\in\mcB$,
\begin{equation}\label{Dtau}
  D\vartheta(p+x) = D\psi_3(\Phi^{-1}(x))\circ D\Phi^{-1}(x).
\end{equation}

For $x'\in B_{\frac{c_1}{\gamma_p},H_1}(0) $ and  any $y\in H_1$
we have
\begin{equation*}\label{eq:psi3}
  D\psi_3(x')y =
  \bigg(\frac{D\omega(p+x')y}{S(x')}
  -\frac{DS(x')y\,\omega(p+x')}{S(x')^2},
   -\frac{DS(x')y}{S(x')^2}\bigg)^t.
\end{equation*}
Therefore,
\begin{eqnarray*}\label{eq:psi3-norm}
\|D\psi_3(x')\| &\leq &
\frac{\|D\omega(p+x')\|}{S(x')}+ \frac{\|DS(x')\|\|\omega(p+x')\|}{S(x')^2}
+\frac{\|DS(x')\|}{S(x')^2}\notag\\
&\underset{S(x')\geq1}{\leq}&
{2.3}\|x'\|\gamma_p +
1.0016\|x'\|{1.15}\|x'\|^2\gamma_p
+ 1.0016\|x'\|\\
&\underset{\|x'\|\le \frac{c_1}{\gamma_p}, \gamma_p\geq1}{\leq}& \|x'\|\gamma_p \big({2.3}+1.0016\cdot {1.15}\cdot c_1^2
+1.0016\big)\;\leq\; {3.303}\|x'\|\gamma_p \notag
\end{eqnarray*}
by Proposition~\ref{prop:M2} and Inequality~\eqref{DS}.

Going back to \eqref{Dtau},
using that $D\Phi^{-1}(x) =
(D\Phi(\Phi^{-1}(x)))^{-1}$, the above inequality  and
Lemma~\ref{lem:inject}(i,ii), we obtain for any $x\in \mcB$,
\begin{eqnarray*}
  \|D\vartheta(p+x)\| &\leq &\|D\psi_3(\Phi^{-1}(x))\|\,\|D\Phi^{-1}(x)\| \\
&\leq&
  \|D\psi_3(\Phi^{-1}(x))\| \big\|(D\Phi(\Phi^{-1}(x)))^{-1}\big\|
 \  \leq \  {3.303}\|\Phi^{-1}(x)\|\gamma_p\cdot 1.0013 \\ & \le &
  {3.303}\cdot 1.0004 \|x\|\gamma_p\cdot 1.0013
 \  \leq \ {3.4}\|x\|\gamma_p.
\end{eqnarray*}

Now, we deduce that
$$
   \|\vartheta(p+x)\|\leq {1.7}\|x\|^2\gamma_p.
$$
the same way we deduced the bound for $\|\omega(p+x)\|$
in Proposition~\ref{prop:M2}.

Finally, Lemma~\ref{lem:inject} also implies that
$B_{\frac{c_2}{\gamma_p},H_1}(0) \subset
\Phi\Big(B_{\frac{c_1}{\gamma_p},H_1}(0)\Big)$,
since for $\|x'\|=\frac{c_1}{\gamma_p}$,
$\|\Phi(x')\|\ge \frac{0.9997c_1}{\gamma_p}
\ge \frac{c_2}{\gamma_p}$.
\eproof

\begin{lemma}\label{lem:M4}
Let $\varphi:H_1\to H_3$ be any linear map and
$E\subset H_1\times H_3$ be the graph of $\varphi$.
Then,
\begin{description}
\item[(i)] $E^\perp=\{(-\varphi^*(v),v)\mid v\in H_3\}
\subset H_1\times H_3$.
\item[(ii)] Let $
   w\in H_3\cap \Big((p+x,\vartheta(p+x))+E^\perp\Big)$
for $\vartheta$ the map of Proposition~\ref{prop:M3} and
$x\in \Phi\Big(B_{\frac{c_1}{\gamma_p},H_1}(0)\Big)\subset H_1$.
Then  $\|w-p\|\geq \frac{\|x\|}{\|\varphi\|}-\|\vartheta(p+x)\|$.
\end{description}
\end{lemma}

\proof
(i) For all $x\in H_1$ and $v\in H_3$ we have
$$
  \langle (x,\varphi(x)),(-\varphi^*(v),v)\rangle
 =\langle x,-\varphi^*(v)\rangle+\langle\varphi(x),v\rangle
 =-\langle x,\varphi^*(v)\rangle+\langle x,\varphi^*(v)\rangle
 =0.
$$
This shows that the linear space
$\{(-\varphi^*(v),v)\mid v\in H_3\}$, of dimension $\dim(H_3)$,
is included in $E^\perp$. The reverse inclusion follows as both
spaces have the same dimension.
\medskip

\noindent (ii)
As
$w\in \big((p+x,\vartheta(p+x))+E^\perp\big)
=\big((x,p+\vartheta(p+x))+E^\perp\big)$, we use
Lemma~\ref{lem:M4}(i) to deduce the existence of
$v\in H_3$ such that
$w=(x,p+\vartheta(p+x))+(-\varphi^*(v),v)\in H_1\times H_3$.
Hence, since $w\in H_3$,
$x-\varphi^*(v)=0$, i.e., $x=\varphi^*(x)$,  and
$w=p+\vartheta(p+x)+v$, i.e.,
$w-p=\vartheta(p+x)+v$. We deduce
$$
     \|v\|\geq \frac{\|x\|}{\|\varphi^*\|} = \frac{\|x\|}{\|\varphi\|}
$$
and, consequently,
$\|w-p\|\geq \|v\|-\|\vartheta(p+x)\| \geq
\frac{\|x\|}{\|\varphi\|}-\|\vartheta(p+x)\|$.
\eproof

\proofof{Theorem~\ref{prop:tau-gamma}}
We show that for all points $p,q\in\mcM_\IS$ the
normals $N_p$ and $N_q$ of $\mcM_\IS$ at $p$ and $q$, i.e.,
the normal spaces to  their tangent planes at $\mcM_\IS$,
either do not intersect or, if they do, the intersection points
lie outside
$B_{\frac{c_2}{2\gamma_p}}(p) \cap B_{\frac{c_2}{2\gamma_p}}(q)$.
Therefore,
$$
  {\tau(f)}\ge \min_{p\in \mcM_\IS}  \frac{c_2}{2\, \gamma_p}
  =  \frac{c_2}{2\, \max_{p\in\mcM_\IS}\gamma_p} \ge
   \frac{1}{87\, \Gamma(f)},
$$
since $\mcM_\IS$ is compact and $c_2=0.023$.

To prove this statement, we take $p$ to be the point in the
preceding development (which is arbitrary on $\mcM_\IS$)
 and divide by cases.
\smallskip

\noindent (i) If $\|q-p\|\ge  \frac{c_2}{\gamma_p}$, then
$B_{\frac{c_2}{2\gamma_p}}(p)\cap
B_{\frac{c_2}{2\gamma_p}}(q)=\emptyset$, which implies that
the normals $N_p$ and $N_q$ cannot intersect at any point in
the intersection of these two balls.

\noindent (ii)
If $\|q-p\|<\frac{c_2}{\gamma_p}$, then
$q\in \mcM_\IS \cap \phi\Big(B_{\frac{c_2}{\gamma_p}}(p)\Big) $
is in the hypothesis of Proposition~\ref{prop:M3}. Let
$x_0\in \Phi\Big(B_{\frac{c_1}{\gamma_p},H_1}(0) \Big) \subset H_1$
be such that $q=(p+x_0,\vartheta(p+x_0))$. Then
$\Big(\frac{c_2}{\gamma_p}\Big)^2>\|q-p\|^2
= \|x_0\|^2+ \|\vartheta(p+x_0)\|^2 \ge \|x_0\|^2$ implies
$x_0 \in B_{\frac{c_2}{\gamma_p},H_1}(0)$, and hence, by the last
statement in Proposition~\ref{prop:M3}, $p+x_0$ belongs
to the domain of $\vartheta$ and we may consider its derivative
$$
   \varphi:=D\vartheta(p+x_0): H_1\to H_3.
$$

Then the graph $E:=\Gr(\varphi)$ is a linear subspace of
$\R^{n+1}$ and the normal
$N_q$ to $E$ at $q=(p+x_0,\vartheta(p+x_0))$ equals
$(p+x_0,\vartheta(p+x_0))+E^\perp$. Analogously the normal
$N_p$ of $\mcM_\IS$ at $p$ equals $p + H_1^{\perp} = p+H_3=H_3$.

Suppose now that $N_q=(p+x_0,\vartheta(p+x_0))+E^\perp$ intersects
$N_p=H_3$ at a point $w$.
Applying Lemma~\ref{lem:M4}(ii) and
Proposition~\ref{prop:M3} we obtain
\begin{eqnarray*}
 \|w-p\| &\geq&
         \frac{\|x_0\|}{\|D\vartheta(p+x_0)\|}-\|\vartheta(p+x_0)\|
  \,\geq\, \frac{\|x_0\|}{3.4\,\gamma_p\|x_0\|}-1.7\,\gamma_p\|x_0\|^2 \\
    &\geq & \frac{1}{3.4\,\gamma_p}-\frac{1.7\,c_2^2}{\gamma_p}
    \;=\;\frac{c_2}{2\gamma_p}\Big(\frac1{1.7\,c_2}-3.4\,c_2\Big)
    \,\geq\, \frac{c_2}{2\gamma_p}
\end{eqnarray*}
the third inequality as $\|x_0\|\leq \frac{c_2}{\gamma_p}$.
This shows that $N_p$ and $N_q$ do not  intersect
in $B_\frac{c_2}{2\gamma_p}(p)$.
\eproof

\section{On numerical stability}\label{sec:stability}

In this last section we deal with the numerical stability of our
algorithms. Part~(iv) of Theorem~\ref{thm:Main} claims that our
algorithms are numerically stable. We now give  a
precise meaning to this claim.

Numerical stability refers to the effects of finite-precision
arithmetic in the final result of a computation. During the execution
of such computation real numbers $x$ are systematically replaced by
approximations $\fl(x)$ satisfying that
$$
     \fl(x)=x(1+\delta), \qquad\mbox{with $|\delta|\leq\emac$}
$$
where $\emac\in(0,1)$ is the {\em machine precision}. If the algorithm
is computing a function $\varphi:\R^p\to\R^q$ a common
definition of stability says that the algorithm is {\em forward
  stable} when, for sufficiently small $\emac$ and for each input
$a\in\R^p$, the computed point
$\tilde{\varphi(a)}\in\R^q$ satisfies
\begin{equation}\label{eq:forward-stab}
  \Big\|\tilde{\varphi(a)}-\varphi(a)\Big\|\leq \emac\, \|\varphi(a)\| \cond(a) P(p,q).
\end{equation}
Here $P$ is a polynomial (which in practice should be of small degree)
and $\cond(a)$ is the condition number of $a$ given by
\begin{equation}\label{eq:cond}
  \cond(a):=\lim_{\delta\to 0} \sup_{\|\tilde{a}-a\|\leq \delta}
                  \frac{\|\varphi(\tilde{a})-\varphi(a)\|}{\|\tilde{a}-a\|}\,
                  \frac{\|a\|}{\|\varphi(a)\|}.
\end{equation}
We observe that $\cond(a)$ depends on $\varphi$ and $a$
but not on the
algorithm and that inequality~\eqref{eq:forward-stab} is
satisfied in first order whenever the algorithm is {\em backward stable},
that is, whenever it satisfies that
\begin{equation}\label{eq:backward-stab}
  \tilde{\varphi(a)}=\varphi(\tilde{a}), \qquad
  \mbox{for some $\tilde{a}$ satisfying $\|\tilde{a}-a\|\leq \|a\|\emac\, P(p,q)$}.
\end{equation}

These notions are appropriate for a continuous function $\varphi$
(such as in matrix inversion, the solution of linear systems of
equations, the computation of eigenvalues, \dots) but
 not  for discrete-valued problems:  if the range of $\varphi$
is discrete (as in deciding the feasibility of a linear program,
counting the number of solutions of a polynomial system, or computing
Betti numbers) then definition~\eqref{eq:cond} becomes meaningless
(see~\cite[Overture, \S6.1, and \S9.5]{Condition} for a detailed exposition of these issues). For these,
a now common definition of condition number, pioneered by Jim
Renegar~\cite{Renegar94b,Renegar95,Renegar95b}, consists of
identifying the set $\Sigma$ of ill-posed inputs and taking the
condition of $a$ as
the relativized inverse of the distance from $a$ to $\Sigma$. That is,
one takes
\begin{equation}\label{eq:cond-dec}
  \msC(a):=\frac{\|a\|}{\dist(a,\Sigma)}.
\end{equation}
Proposition~\ref{prop:kappa-dist} shows that our condition
number $\kappa(f)$
is bounded by such an expression (with respect to the set of ill-posed
inputs $\Sigma_{\R}$). 

 The idea of
stability changes together with the definition of condition. The issue now is not the one
underlying~\eqref{eq:forward-stab} ---given $\emac$, how good is
the computed value--- but a different one: {\em how small does $\emac$ need
to be to ensure that the computed output is correct?} The answer to
this question depends on the condition of the input at hand, a
quantity that is generally not known a priori, and
stability results can be broadly divided in two classes. In a
{\em fixed-precision} analysis the algorithm runs with a
pre-established machine precision and the users have no guarantee that
the returned output is correct. They only know that if the input $a$ is
well conditioned (i.e., smaller than a bound depending on $\emac$)
then the answer is correct. In a {\em variable-precision} analysis the
algorithm has the capacity to adjust its machine precision during the
execution and returns an output which is guaranteed to be
correct. Needless to say, not all algorithms may be brought to a
variable-precision analysis. But in the last decades a number of
problems such as feasibility for semialgebraic systems~\cite{CS98}
or for linear programs~\cite{CP01}, real zero counting of polynomial
systems~\cite{CKMW1}, or the computation of optimal bases for
linear programs~\cite{ChC03} have been given such analysis.

In all these cases, it is shown that the finest precision $\emac^*$
used by the algorithm satisfies
\begin{equation}\label{eq:stab-dec}
  \emac^*=\frac{1}{(p \msC(a))^{\Oh(1)}}
\end{equation}
where $p$ is the {\em size} of the input and $\msC(a)$
is the condition number defined in~\eqref{eq:cond-dec}. We can
(and will) consider algorithms satisfying~\eqref{eq:stab-dec} to be
{\em stable} as this bound implies that the number of bits in the
mantissa of the floating-point numbers occurring in the computation
with input $a\in\R^p$ is bounded by $\Oh(\log_2 p+\log_2  \msC(a))$.

It is in this sense that our algorithms are stable.

\begin{proposition}
The algorithms in Propositions~\ref{prop:betti-S}
and~\ref{prop:betti-P} computing the homology groups of
spherical and projective sets, respectively, can be modified
to work with variable-precision and satisfy the following.
Their cost, for an input $f\in\Hd[m]$, remain
$$
 (n D \kappa(f))^{\Oh(n^2)}
$$
and the finest precision $\emac^*$ used by the algorithm
is
$$
  \emac^*=\frac{1}{(n D \kappa(f)\log N)^{\Oh(1)}}.
$$
\end{proposition}

\sketchproof
A key observation for the needed modification is that only
the routine {\sf Covering} needs to work with
finite precision. Indeed, we can modify this routine
to return a pair $\{\mcX,\e\}$ where all numbers, coordinates of
points $x$ in $\mcX$ and $\e$, are rational numbers
(expressed as quotients of integers in binary form).
Furthermore, we can do so such that the differences $\|x-\tilde{x}$
and $|\e-\tilde\e|$ between the real objects and their rational
approximations are small. Sufficiently small actually for
Proposition~\ref{prop:SNW} to apply to $(\tilde\mcX,\tilde\e)$
(recall that Remark~\ref{rem:tau} gives us plenty of room to do so).

From this point on, the computation of the nerve $\mcN$  and then of
the homology groups of either $\mcM_{\IS}$ or $\mcM_{\P}$  is done
symbolically (i.e.,  with infinite precision). The complexity of the whole
procedure, that is, its cost, which now takes account of the
size of the rational numbers occuring during the computation,
remains within the same general bound
in the statement.

We therefore only need to show that a variable-precision version
of {\sf Covering} can be devised that returns an output with rational
components and that satisfies the bounds in the statement. This
version is constructed, essentially, as the variable-precision version
of the algorithm for counting roots in~\S5.2 of~\cite{CKMW1} is
constructed in~\S6.3 of that paper. We do not give all the details here
since these  do not add anything new to our understanding
of the algorithm: we just    ``make room'' for errors
by weakening the desired inequalities by a factor of~2; in our case,
the inner loop of the algorithm becomes
\bigskip

\codigo
\bodycode{
\espacio
for all $x\in\mcG_\eta$ \\
\eespacio if $\ba(f,x)\leq\frac{\a_0}{2}$ and
$\frac{1}{1000\,\bg(f,x)}\geq r$
   and $4.4\,\bb(f,x)< r$ then\\
\eeespacio $\mcX:=\mcX\cup \{x\}$\\
\eespacio elsif $\|f(x)\|\geq 2\delta(f,\eta)$ then do nothing\\
\eespacio elsif go to (*)\\
\espacio return the pair $\{\mcX,\e\}$ and halt\\
\espacio end for
}
\fcodigo

\noindent
Also, as Proposition~\ref{prop:SNW} does neither require
the points of $\mcX$ to belong to the sphere, nor a precise
value for $\e$, there is no harm in returning points (with
rational coefficients) close to the sphere and to work with
a good (rational) approximation $\e$ of $3.5\sqrt{\sep(\eta)}$.
\eproof

We close this section by recalling that the biggest mantissa
required in a floating-point computation with input $f$ has
$\Oh(\log_2 (n D \kappa(f)\log N))$ bits.
If $f$ is randomly drawn from $\IS^{N-1}$ this is a random
variable. Using the second bound in
Theorem~\ref{th:BCL} along with
Propositions~\ref{prop:kappa-dist} and~\ref{prop:discrim}
it follows that the expectation for the number of bits in
this longest mantissa is of the order of
$$
     \Oh\big(n\log_2(Dm) +\log_2 N +\log_2 n\big).
$$
This is a relatively small quantity compared with
(and certainly polynomially bounded in) the size $N$ of input $f$.

{\small

\begin{thebibliography}{10}

\bibitem{AllGeo90}
E.L. Allgower and K.~Georg.
\newblock {\em Numerical Continuation Methods}.
\newblock Springer-Verlag, 1990.

\bibitem{wacco}
D.~Amelunxen and M.~Lotz.
\newblock Average-case complexity without the black swans.
\newblock To appear at {\em \JoC}. Available at {\tt arXiv:1512.09290}, 2016.

\bibitem{Basu:08}
S.~Basu.
\newblock Computing the top {B}etti numbers of semialgebraic sets defined by
  quadratic inequalities in polynomial time.
\newblock {\em Found. Comput. Math.}, 8(1):45--80, 2008.

\bibitem{BaPoRo:08}
S.~Basu, R.~Pollack, and M.-F. Roy.
\newblock Computing the first {B}etti number of a semi-algebraic set.
\newblock {\em Found. Comput. Math.}, 8(1):97--136, 2008.

\bibitem{Bjo:95}
A.~Bj\"orner.
\newblock Topological methods.
\newblock In R.~Graham, M.~Grotschel, and L.~Lovasz, editors, {\em Handbook of
  Combinatorics}, pages 1819--1872. North-Holland, Amsterdam, 1995.

\bibitem{BCSS98}
L.~Blum, F.~Cucker, M.~Shub, and S.~Smale.
\newblock {\em Complexity and Real Computation}.
\newblock Springer-Verlag, 1998.

\bibitem{BSS89}
L.~Blum, M.~Shub, and S.~Smale.
\newblock On a theory of computation and complexity over the real numbers:
  {NP}-completeness, recursive functions and universal machines.
\newblock {\em \BAMS}, 21:1--46, 1989.

\bibitem{BC03}
P.~B\"urgisser and F.~Cucker.
\newblock Counting complexity classes for numeric computations {II}: Algebraic
  and semialgebraic sets.
\newblock {\em \JoC}, 22:147--191, 2006.

\bibitem{BC09}
P.~B\"urgisser and F.~Cucker.
\newblock Exotic quantifiers, complexity classes, and complete problems.
\newblock {\em Found. Comput. Math.}, 9:135--170, 2009.

\bibitem{Condition}
P.~B\"urgisser and F.~Cucker.
\newblock {\em Condition}, volume 349 of {\em Grundlehren der mathematischen
  Wissenschaften}.
\newblock Springer-Verlag, Berlin, 2013.

\bibitem{ChC03}
D.~Cheung and F.~Cucker.
\newblock Solving linear programs with finite precision: {II}. {A}lgorithms.
\newblock {\em \JoC}, 22:305--335, 2006.

\bibitem{Collins}
G.E. Collins.
\newblock {\em Quantifier elimination for real closed fields by cylindrical
  algebraic deccomposition}, volume~33 of {\em \LNCS}, pages 134--183.
\newblock Springer-Verlag, 1975.

\bibitem{Cucker99b}
F.~Cucker.
\newblock Approximate zeros and condition numbers.
\newblock {\em \JoC}, 15:214--226, 1999.

\bibitem{CDW:05}
F.~Cucker, H.~Diao, and Y.~Wei.
\newblock Smoothed analysis of some condition numbers.
\newblock {\em Numer. Lin. Alg. Appl.}, 13:71--84, 2006.

\bibitem{CKMW1}
F.~Cucker, T.~Krick, G.~Malajovich, and M.~Wschebor.
\newblock A numerical algorithm for zero counting. {I}: {C}omplexity and
  accuracy.
\newblock {\em \JoC}, 24:582--605, 2008.

\bibitem{CKMW2}
F.~Cucker, T.~Krick, G.~Malajovich, and M.~Wschebor.
\newblock A numerical algorithm for zero counting. {II}: {D}istance to
  ill-posedness and smoothed analysis.
\newblock {\em J. Fixed Point Theory Appl.}, 6:285--294, 2009.

\bibitem{CKMW3}
F.~Cucker, T.~Krick, G.~Malajovich, and M.~Wschebor.
\newblock A numerical algorithm for zero counting. {III}: {R}andomization and
  condition.
\newblock {\em Adv. Applied Math.}, 48:215--248, 2012.

\bibitem{CP01}
F.~Cucker and J.~Pe\~na.
\newblock A primal-dual algorithm for solving polyhedral conic systems with a
  finite-precision machine.
\newblock {\em \SIOPT}, 12:522--554, 2002.

\bibitem{CS98}
F.~Cucker and S.~Smale.
\newblock Complexity estimates depending on condition and round-off error.
\newblock {\em \JACM}, 46:113--184, 1999.

\bibitem{DKS:13}
Carlos D'Andrea, Teresa Krick, and Mart{\'{\i}}n Sombra.
\newblock Heights of varieties in multiprojective spaces and arithmetic
  {N}ullstellens\"atze.
\newblock {\em Ann. Sci. \'Ec. Norm. Sup\'er. (4)}, 46(4):549--627 (2013),
  2013.

\bibitem{Demmel88}
J.~Demmel.
\newblock The probability that a numerical analysis problem is difficult.
\newblock {\em Math. Comp.}, 50:449--480, 1988.

\bibitem{EdelHar:10}
H.~Edelsbrunner and J.L. Harer.
\newblock {\em Computational topology}.
\newblock American Mathematical Society, Providence, RI, 2010.
\newblock An introduction.

\bibitem{Kostlan88}
E.~Kostlan.
\newblock Complexity theory of numerical linear algebra.
\newblock {\em J. of Computational and Applied Mathematics}, 22:219--230, 1988.

\bibitem{Lotz15}
M.~Lotz.
\newblock On the volume of tubular neighborhoods of real algebraic varieties.
\newblock {\em Proc. Amer. Math. Soc.}, 143(5):1875--1889, 2015.

\bibitem{NiSmWe08}
P.~Niyogi, S.~Smale, and S.~Weinberger.
\newblock Finding the homology of submanifolds with high confidence from random
  samples.
\newblock {\em Discrete Comput. Geom.}, 39:419--441, 2008.

\bibitem{NoTo:15}
V.~Noferini and A.~Townsend.
\newblock Numerical instability of resultant methods for multidimensional
  rootfinding.
\newblock To appear at {\em SIAM J. Num. Analysis}. Available at {\tt
  arXiv:1507.00272}.

\bibitem{Ren92a}
J.~Renegar.
\newblock On the computational complexity and geometry of the first-order
  theory of the reals. {P}art {I}.
\newblock {\em Journal of Symbolic Computation}, 13:255--299, 1992.

\bibitem{Renegar94b}
J.~Renegar.
\newblock Some perturbation theory for linear programming.
\newblock {\em \MP}, 65:73--91, 1994.

\bibitem{Renegar95}
J.~Renegar.
\newblock Incorporating condition measures into the complexity theory of linear
  programming.
\newblock {\em \SIOPT}, 5:506--524, 1995.

\bibitem{Renegar95b}
J.~Renegar.
\newblock Linear programming, complexity theory and elementary functional
  analysis.
\newblock {\em \MP}, 70:279--351, 1995.

\bibitem{Scheib:07}
P.~Scheiblechner.
\newblock On the complexity of deciding connectedness and computing {B}etti
  numbers of a complex algebraic variety.
\newblock {\em J. Complexity}, 23(3):359--379, 2007.

\bibitem{Scheib:12}
P.~Scheiblechner.
\newblock Castelnuovo-{M}umford regularity and computing the de {R}ham
  cohomology of smooth projective varieties.
\newblock {\em Found. Comput. Math.}, 12(5):541--571, 2012.

\bibitem{Bez1}
M.~Shub and S.~Smale.
\newblock Complexity of {B}\'ezout's {T}heorem {I}: geometric aspects.
\newblock {\em \JAMS}, 6:459--501, 1993.

\bibitem{Bez2}
M.~Shub and S.~Smale.
\newblock Complexity of {B}\'ezout's {T}heorem {II}: volumes and probabilities.
\newblock In F.~Eyssette and A.~Galligo, editors, {\em Computational Algebraic
  Geometry}, volume 109 of {\em Progress in Mathematics}, pages 267--285.
  Birkh\"auser, 1993.

\bibitem{Bez3}
M.~Shub and S.~Smale.
\newblock Complexity of {B}\'ezout's {T}heorem {III}: condition number and
  packing.
\newblock {\em Journal of Complexity}, 9:4--14, 1993.

\bibitem{Bez5}
M.~Shub and S.~Smale.
\newblock Complexity of {B}\'ezout's {T}heorem {V}: polynomial time.
\newblock {\em \TCS}, 133:141--164, 1994.

\bibitem{Bez4}
M.~Shub and S.~Smale.
\newblock Complexity of {B}\'ezout's {T}heorem {IV}: probability of success;
  extensions.
\newblock {\em SIAM J. of Numer. Anal.}, 33:128--148, 1996.

\bibitem{Smale86}
S.~Smale.
\newblock Newton's method estimates from data at one point.
\newblock In R.~Ewing, K.~Gross, and C.~Martin, editors, {\em The Merging of
  Disciplines: New Directions in Pure, Applied, and Computational Mathematics}.
  Springer-Verlag, 1986.

\bibitem{stor:96}
A.~Storjohann.
\newblock Nearly optimal algorithms for computing {S}mith normal forms of
  integer matrices.
\newblock In {\em Proceedings of the {I}nternational {S}ymposium on {S}ymbolic
  and {A}lgebraic {C}omputation ({ISSAC}'96)}, pages 267--274. ACM Press, 1996.

\bibitem{Wut76}
H.R. W\"uthrich.
\newblock Ein {E}ntscheidungsverfahren f\"ur die {T}heorie der
  reell-\-ab\-ge\-schlos\-se\-nen {K}\"orper.
\newblock In E.~Specker and V.~Strassen, editors, {\em Komplexit\"at von
  Entscheidungsproblemen}, volume~43 of {\em \LNCS}, pages 138--162.
  Springer-Verlag, 1976.

\end{thebibliography}

}

\end{document}